\documentclass{article}
\usepackage{amssymb, amsthm, latexsym, amscd, amsfonts, verbatim, amsmath}
\usepackage[all, cmtip]{xy}

\DeclareMathOperator{\cor}{cor}
\DeclareMathOperator{\img}{im}
\DeclareMathOperator{\im}{im}

\DeclareMathOperator{\spec}{Spec}

\DeclareMathOperator{\coker}{coker}
\DeclareMathOperator{\cok}{coker}

\DeclareMathOperator{\Hom}{Hom}

\DeclareMathOperator{\Hilb}{Hilb}
\DeclareMathOperator{\Hilbe}{Hilb_{exc}}

\DeclareMathOperator{\reg}{reg}
\DeclareMathOperator{\DT}{DT}
\DeclareMathOperator{\PT}{PT}
\DeclareMathOperator{\TP}{\pi-PT}
\DeclareMathOperator{\coh}{coh}
\DeclareMathOperator{\supp}{supp}
\DeclareMathOperator{\SSS}{SS}

\DeclareMathOperator{\Hilba}{Hilb^{\alpha}}
\DeclareMathOperator{\St}{St}

\DeclareMathOperator{\HHHH}{H}
\DeclareMathOperator{\Ad}{Ad}

\newcommand{\iso}{\cong}

\newcommand{\wt}[1]{\widetilde{#1}}
\newcommand{\pr}{^{\prime}}
\newcommand{\ppr}{^{\prime \prime}}

\newcommand{\too}{\longrightarrow}

\newcommand{\CC}{\mathbb{C}}
\newcommand{\CCC}{\mathcal{C}}

\newcommand{\OO}{\mathcal{O}}
\newcommand{\HHH}{\mathcal{H}}
\newcommand{\MM}{\mathcal{M}}
\newcommand{\NN}{\mathcal{N}}
\newcommand{\NNN}{\mathbb{N}}
\newcommand{\UU}{\mathcal{U}}
\newcommand{\HH}{\mathcal{H}}

\newcommand{\XX}{\mathcal{X}}

\newcommand{\ZZ}{\mathbb{Z}}

\newcommand{\LL}{\mathbb{L}}
\newcommand{\QQ}{\mathbb{Q}}
\newcommand{\QQQ}{\mathcal{Q}}
\newcommand{\PP}{\mathbb{P}}
\newcommand{\PPP}{\mathcal{P}}

\newcommand{\AAA}{\mathbb{A}}
\newcommand{\AAAA}{\mathcal{A}}

\newcommand{\sub}{\subset}

\newcommand{\Qp}{\mathcal{Q}_{\pi}}
\newcommand{\Pp}{\mathcal{P}_{\pi}}

\newcommand{\mup}{\mu_{\pi}}
\newcommand{\HHp}{ \HH^{\pi}}

\newcommand{\HHe}{\HH_{exc}}
\newcommand{\Hleq}{\HHH}

\newcommand{\oneQp}{1_{\Qp}}
\newcommand{\onePp}{1_{\Pp}}
\newcommand{\oneQpO}{1^{\OO}_{\Qp}}
\newcommand{\onePpO}{1^{\OO}_{\Pp}}

\newcommand{\onePponeO}{1^{\OO}_{\Pp^1}}

\newcommand{\one}[1]{1_{#1}}
\newcommand{\oneo}[1]{1^{\OO}_{#1}}
\newcommand{\oo}{\infty}
\newcommand{\twoo}{\frac{\oo}{2}}
\newcommand{\smloo}{\SSS ( \mu \leq \square \leq \oo)}
\newcommand{\smltwoo}{\SSS ( \mu \leq \square < \twoo)}

\newcommand{\HA}{\HHHH (\CCC)}
\newcommand{\HRA}{\HHHH_{\reg} (\CCC)}

\newcommand{\HC}{\HHHH (\CCC)}
\newcommand{\HRC}{\HHHH_{\reg} (\CCC)}
\newcommand{\HSC}{\HHHH_{\hbox{\tiny sc}} (\CCC)}
\newcommand{\ps}{$\pi$-stable}
\newcommand{\Hilbpbn}{\operatorname{\pi -Hilb}^{(\beta,n)}}
\newcommand{\Hilbp}{\operatorname{\pi-Hilb}}
\newcommand{\Hilbpa}{\operatorname{\pi-Hilb}^{\alpha}}
\newcommand{\Kvar}{K(\text{Var}/\CC)}
\newcommand{\KSt}{K(\text{St}/\CC)}

\newtheorem{example}[equation]{Example}

\newtheorem{lemma}[equation]{Lemma}

\newtheorem{definition}[equation]{Definition}

\newtheorem{prop}[equation]{Proposition}

\newtheorem{theorem}[equation]{Theorem}
\newtheorem{corollary}[equation]{Corollary}
\newtheorem{conjecture}[equation]{Conjecture}

\theoremstyle{definition} \newtheorem{remark}[equation]{Remark}

\newenvironment{pf}{\noindent {\bf Proof:}}{  \hfill $\blacksquare$ \\}
\newenvironment{fpf}{\noindent {\bf Fake proof:}}{  \hfill $\sqcup$ \\}

\title{Curve counting invariants for crepant resolutions}
\author{Jim Bryan and David Steinberg}

\begin{document}
\maketitle

\abstract{We construct curve counting invariants for a Calabi-Yau
threefold $Y$ equipped with a dominant birational morphism $\pi:Y \to
X$. Our invariants generalize the stable pair invariants of
Pandharipande and Thomas which occur for the case when $\pi:Y\to Y$ is
the identity. Our main result is a PT/DT-type formula relating the
partition function of our invariants to the Donaldson-Thomas partition
function in the case when $Y$ is a crepant resolution of $X$, the
coarse space of a Calabi-Yau orbifold $\mathcal{X}$ satisfying the
hard Lefschetz condition. In this case, our partition function is
equal to the Pandharipande-Thomas partition function of the orbifold
$\mathcal{X}$.} Our methods include defining a new notion of stability
for sheaves which depends on the morphism $\pi $. Our notion
generalizes slope stability which is recovered in the case where $\pi
$ is the identity on $Y$.

\tableofcontents


\section{Introduction} Donaldson-Thomas (DT) theory of a Calabi-Yau
threefold $X$ gives rise to subtle deformation invariants. They are
considered to be the mathematical counterparts of BPS state counts in
topological string theory compactified on $X$. Principles of physics
(see \cite{vaf}, \cite{zas}) indicate that the string theory of an
orbifold Calabi-Yau threefold and that of its crepant resolution ought
to be equivalent, so one expects that the DT theories of an orbifold
and its crepant resolution to be equivalent in some way.  In the case
where the orbifold satisfies the hard Lefschetz condition, the crepant
resolution conjecture of \cite{bcy} gives a formula determining the DT
invariants of the orbifold in terms of the DT invariants of the
crepant resolution.

In this article, we begin a program to prove the crepant resolution
conjecture using Hall algebra techniques inspired by those of
Bridgeland \cite{bri2}. In the process, we construct curve counting
invariants for a Calabi-Yau threefold $Y$ equipped with a birational
morphism $\pi:Y \to X$. Our invariants generalize the stable pair
invariants of Pandharipande and Thomas which occur for the case when
$\pi:Y\to Y$ is the identity. Our main result is a PT/DT-type formula
relating the partition function of our invariants to the
Donaldson-Thomas partition function in the case when $Y$ is a crepant
resolution of $X$, the coarse space of a Calabi-Yau orbifold
$\mathcal{X}$ satisfying the hard Lefschetz condition. In this case,
our partition function is equal to the Pandharipande-Thomas partition
function of the orbifold $\mathcal{X}$.


\subsection*{Donaldson-Thomas theory}
Let $Y$ be a smooth projective Calabi-Yau threefold.  Let $K(Y)$ be the numerical K-theory of $Y$, i.e. the quotient of the K-group of $\coh (Y)$ by the kernel of the Chern character map to cohomology. The Hilbert scheme of $Y$, $\Hilba (Y)$, parametrizes quotients $\OO_Y \to \OO_Z$, such that the class of $\OO_Z$ in $K(Y)$ is $\alpha$. The group $K(Y)$ is filtered by the dimension of the support:
$$
F_0 K(Y) \sub F_1 K(Y) \sub F_2 K(Y) \sub F_3 K(Y) = K(Y).
$$
In this article, we will focus on curves, i.e., $\alpha\in F_1K(Y)$, with $\text{ch}(\alpha) = (0,0,\beta, n)$, where $\beta\in H^4(Y, \ZZ)$ is a curve class, and $n\in H^6(Y, \ZZ)\iso \ZZ$ is the holomorphic Euler characteristic. In \cite{tho}, an obstruction theory for this moduli space is constructed, which produces (by \cite{bf1}) a virtual fundamental cycle. Donaldson-Thomas invariants are defined by integrating over the zero-dimensional virtual fundamental class:

$$
DT^\alpha (Y)= \int_{[\Hilba (Y)]^{vir}} 1. 
$$

Since the obstruction theory is symmetric, we may also express the invariants as the Euler characteristic of $\Hilba (Y)$ weighted by Behrend's microlocal function \cite{beh}:
$$
DT^\alpha (Y) = \sum_{n\in \ZZ} n \chi(\nu^{-1}(n)),
$$
where $\nu: \Hilba (Y) \to \ZZ$ is Behrend's function.

Following \cite{mnop}, we assemble the invariants into a partition function 
$$
\DT (Y) = \sum_{\alpha \in F_1K(Y)} \DT^\alpha (Y) q^\alpha.
$$

\begin{remark} In \cite{joso}, Donaldson-Thomas invariants are greatly generalized, from the case of structure sheaves of curves to that of arbitrary sheaves. The price of admission to this generality is the formidable machinery of Joyce \cite{joy1, joy2, joy3, joy4, joy5}. An even more ambitious program of generalization is being lead by Kontsevich and Soibelman \cite{koso}.
\end{remark}
\subsection*{The Donaldson-Thomas crepant resolution conjecture}

We follow \cite{bcy} in our treatment of the crepant resolution conjecture.

An \emph{orbifold CY3} is defined to be a smooth, quasi-projective, Deligne-Mumford stack $\XX$ over $\CC$ of dimension three having generically trivial stabilizers and trivial canonical bundle, 
$$
K_\XX \iso \OO_\XX.
$$
The definition implies that the local model for $\XX$ at a point $p$
is $[\CC^3 / G_p]$ where $G_p\sub SL(3, \CC)$ is the (finite) group of
automorphisms of $p$. The orbifold CY3s that appear in this article
will all be projective and satisfy the \emph{hard Lefschetz condition}
\cite[definition 1.1]{brgr}, which in this case is equivalent
\cite[lemma 24]{brgh} to the condition that all $G_p$ are finite
subgroups of $SO(3) \sub SU(3)$ or $SU(2)\sub SU(3)$.

Let $X$ denote the coarse space of $\XX$. A \emph{crepant resolution} of $X$ is a resolution of singularities $\pi:Y\to X$ such that $\pi^*K_X \iso K_Y$. Lemma (1) and Proposition (1) of \cite{vie} prove that 
\begin{equation} \label{rationalsing}
R^{\bullet}\pi_*\OO_Y \cong  \OO_X.
\end{equation} 
The results of \cite{bkr} and \cite{chts} prove that one distinguished crepant resolution of $X$ is 
\begin{equation} 
\label{res}
Y = \Hilb ^{[\OO_p]}(\XX),
\end{equation}
the Hilbert scheme parametrizing substacks in the class $[\OO_p] \in
F_0K(\XX)$. The hard Lefschetz condition implies that the resolution
is \emph{semi-small} (i.e., that the fibres of $\pi$ are zero- or
one-dimensional), and that the singular locus of $X$ is
one-dimensional; see \cite{bosa, brgh}. Furthermore, \cite{bkr} and
\cite{chts} prove that there is a Fourier-Mukai isomorphism
$$
\Psi: D^b(Y) \to D^b(\XX)
$$ 
defined by 
$$
E \mapsto Rq_*p^*E
$$
where
$$
p: Z \to Y, \hbox{    } q:Z \to \XX
$$
are the projections from the universal substack $Z \sub \XX \times Y$ onto each factor. This isomorphism descends to an isomorphism of $K$-theory also denoted $\Psi: K(Y) \to K(\XX)$. It does not respect the filtration by dimension. However, the hard Lefschetz condition implies that the image of $F_0K(\XX)$ is contained in $F_1K(Y)$, under the inverse $\Phi$ of $\Psi$. We call the image $F_{exc}K(Y)$; its elements can be represented by formal differences of sheaves supported on the exceptional fibres of $\pi: Y \to X$. We define the multi-regular part of $K$-theory, $F_{mr}(\XX)$, to be the preimage of $F_1K(Y)$ under $\Psi$. Its elements can be represented by formal differences of sheaves supported in dimension one where at the generic point of each curve in the support, the associated representation of the stabilizer group of that point is a multiple of the regular representation. The following filtrations are respected by $\Psi$:
$$
F_{exc}K(Y) \subset F_1K(Y) \subset K(Y)
$$
$$
F_0K(\XX) \subset F_{mR}K(\XX) \subset K(\XX).
$$

Define the \emph{exceptional} $\DT$ generating series of $Y$, the multi-regular generating series, and degree zero generating series of $\XX$ to be:
$$
\DT_{exc}(Y) = \sum_{\alpha \in F_{exc} K(Y)} DT^{\alpha}(Y)q^\alpha,
$$
$$
\DT_{mr}(\XX) = \sum_{\alpha \in F_{mr}K(\XX)} DT^{\alpha}(\XX)q^\alpha
$$
$$
\DT_0 (\XX) =  \sum_{\alpha \in F_{0}K(\XX)} DT^{\alpha}(\XX)q^\alpha
$$
We state the crepant resolution conjecture of \cite[conjecture 1]{bcy}:
\begin{conjecture} Let $\XX$ be an orbifold CY3 satisfying the hard Lefschetz condition. Let $Y$ be the Calabi-Yau resolution of $\XX$ given by equation~\ref{res}. Then using $\Psi$ to identify the variables, we have an equality
$$
\frac{\DT_{mr}(\XX)}{\DT_0(\XX)} = \frac{\DT(Y)}{\DT_{exc}(Y)}.
$$
\end{conjecture}
\noindent This article makes progress towards proving this conjecture.

In his recent article \cite{cale}, John Calabrese proves a relationship between the DT invariants of a Calabi-Yau threefold and its flop. This problem is similar in many respects to the crepant resolution conjecture studied in this thesis, and Calabrese uses many similar techniques. He constructs a torsion pair and new counting invariants which he relates to invariants on the flop via equations in the Hall algebra and the integration map. While this is very similar to our approach in outline, the actual torsion pair and counting invariants that Calabrese considers (even when adapted to the orbifold setting) are quite different from ours. It would be very interesting to find the precise relationship between the two approaches. An even more recent preprint \cite{cal2} of Calabrese proves the DT crepant resolution conjecture, utilizing his earlier paper \cite{cale}.

\subsection*{$\pi$-stable pairs}
Objects of the Hilbert scheme may be viewed as two-term complexes,
$$
\OO_Y \stackrel{\gamma}{\to} G,
$$
where the cokernel of $\gamma$ must be zero, and where $G$ may be any sheaf admitting such a map $\gamma$. The new invariants introduced in this article, $\pi$-stable pairs, are a modification of this idea. They have been constructed with a view towards proving the crepant resolution conjecture, and as such, they depend on a crepant resolution $Y \stackrel{\pi}{\to} X$ as described in the previous section. The objects of our moduli space allow more variation in our cokernels, but less in the sheaf $G$. In particular, a two-term complex
$$
\OO_Y \stackrel{\gamma}{\to} G
$$
is a $\pi$-\emph{stable pair} (c.f. definition~\ref{psp}) if: 
\begin{enumerate}
\item $R^\bullet\pi_* \cok (\gamma)$ is a zero-dimensional sheaf on $X$, and
\item $G$ admits only the zero map from any sheaf $P$ with the property above, namely that $R^\bullet\pi_*P$ is a zero-dimensional sheaf.
\end{enumerate}

\begin{remark} These pairs were inspired by, and are a generalization of, the stable pairs of Pandharipande and Thomas \cite{pt1}. In fact, when $X=Y$ and $\pi$ is the identity map, the above definition reduces to their definition of stable pairs.
\end{remark}

Below, we prove that there is a finite-type constructible space, $\Hilbpa$ parametrizing these objects with $[G] =\alpha \in K(Y)$. We may then define invariants
$$
\TP^\alpha (Y) = \sum_{n\in \ZZ} n \chi(\nu^{-1}(n)),
$$
where $\nu: \Hilbpa \to \ZZ$ is Behrend's microlocal function. Note that if $\pi: Y \to X$ is the identity, then $\TP^\alpha (Y) =  \PT^\alpha (Y),$ the usual Pandharipande-Thomas invariants of $Y$. As with Donaldson-Thomas theory, we collect the invariants into a generating series, 
$$
\TP (Y) = \sum_{\alpha \in F_1 K(Y)} \TP^\alpha(Y) q^\alpha.
$$

\subsection*{Main result}
The following theorem rests the work of Bridgeland \cite{bri2} and Joyce--Song \cite{joso}, and we therefore require our Calabi-Yau threefold $Y$ to satisfy 
$$
H^1(Y, \OO_Y)=0.
$$ 

\begin{theorem} \label {therecanbeonlyone}
Let $X$ be a projective Calabi-Yau threefold that is the coarse space
of an orbifold CY3 $\XX$ that satisfies the hard Lefschetz
condition. Let $\pi: Y \to X$ be the resolution given by
equation~\ref{res}. Then the generating series for the $\pi$-stable
pair invariants and the DT invariants are related by the equation
$$
\TP(Y) = \frac{\DT(Y)}{\DT_{exc}(Y)}.
$$

\end{theorem}

The aim of this article is to prove this theorem. We summarize the chapters below. In chapter 2, we describe a torsion pair $(\Pp, \Qp$) that is crucial to our definition of $\pi$-stable pairs. We explain the similarities between $\pi$-stable pairs and PT stable pairs and objects of the Hilbert scheme. The chapter ends by establishing results about the moduli space of $\pi$-stable pairs.

In chapter 3, we recall the concept of a stability condition in the sense of Joyce. We then define the stability condition that we will use through out. The rest of the chapter is dedicated to proving that we may apply Joyce's powerful machinery. 

In chapter 4, we introduce the Harder-Narasimhan filtration for our stability condition, which will be our main tool to prove the relationship between the stability condition and the torsion pair from chapter 2.

In chapter 5, we introduce the motivic Hall algebra. 

In chapter 6, we introduce the infinite-type Hall algebra as a purely pedagogical tool. It helps us to give the essence of the idea of many results, without having to concern ourselves with convergence issues, which are handled in the next chapter. 

In chapter 7, we introduce the Laurent Hall algebra, address the convergence issues alluded to in the previous chapter, and prove theorem~\ref{therecanbeonlyone}. 

\begin{remark} To prove the crepant resolution conjecture, we need to prove that $\TP (Y) = \PT (\XX)$ and then use (\cite{ba2}) Bayer's proof of the PT/DT correspondence on $\XX$,
$$
\PT (\XX) = \frac{\DT_{mR}(\XX)}{\DT_0(\XX)}.
$$
The hope is that the Fourier-Mukai isomorphism $\Psi$ takes \ps \hbox{ }pairs (as an object in $D^b(Y)$) to a PT pair on $\XX$.
\end{remark}

\noindent {\bf Acknowledgements.} The authors wishes to thank Arend
Bayer, Brian Conrad, S\'{a}ndor Kov\'{a}cs, Kai Behrend, John
Calabrese, Kalle Karu, Andrew Morrison, and mathoverflow.com for their
helpful conversations. Special thanks are due to Tom Bridgeland for
invaluable comments and suggestions on an early draft of this article,
as well as feedback on a later draft.  We also thank the referee,
whose in depth report helped improve this paper considerably.


\section {$\pi$-stable pairs} In this section, we define $\pi$-stable pairs, and prove some basic results. 

\subsection*{Categorical constructions}
Let $\AAAA$ be an abelian category. Here we recall the notion of torsion pairs. 
\begin{definition}
\label{def:torsion}

Let 
$(\PPP, \QQQ)$ be a pair of full subcategories of $\AAAA$. 
We say $(\PPP, \QQQ)$ is a \textit{torsion pair}
if the following conditions hold.
\begin{itemize}
\item $\Hom(T, F)=0$ for any $T\in \PPP$ and 
$F\in \QQQ$.

\item Any object $E\in \AAAA$
fits into a unique exact sequence, 
\begin{align}
\label{fits}
0 \to T \to E \to F \to 0, 
\end{align}

with $T\in \PPP$ and $F\in \QQQ$.
\end{itemize}
\end{definition}

We borrow the following lemma from Toda \cite{tod}.

\begin{lemma} \label{toda}

Suppose that $\AAAA$ is a noetherian abelian category. 

(i) Let $\PPP \subset \AAAA$ be a full subcategory which 
is closed under extensions and quotients in $\AAAA$. 
Then for $\QQQ=\{E\in \AAAA : \Hom(\PPP, E)=0\}$, the pair 
$(\PPP, \QQQ)$ is a torsion pair on $\AAAA$. 

(ii) Let $\QQQ \subset \AAAA$ be a full subcategory which 
is closed under extensions and subobjects in $\AAAA$. 
Then for $\PPP=\{E\in \AAAA : \Hom(E, \QQQ)=0\}$, the pair 
$(\PPP, \QQQ)$ is a torsion pair on $\AAAA$. 
\end{lemma}
\begin{pf}
We only show (i), as the proof of (ii) is similar. 
Take $E\in \AAAA$ with $E\notin \QQQ$. Then there is 
$T\in \PPP$ and a non-zero morphism 
$T\to E$. Since $\PPP$ is closed under quotients, 
we may assume that $T\to E$ is a monomorphism in $\AAAA$. 
Take an exact sequence in $\AAAA$, 
\begin{align}\label{easylem}
0 \to T \to E \to F \to 0.
\end{align}
By the noetherian property of $\AAAA$ and the assumption
that $\PPP$ is closed under extensions, we may assume that 
there is no $T\subsetneq T' \subset E$ with $T' \in \PPP$. 
Then we have $F\in \QQQ$ and (\ref{easylem}) gives the 
desired sequence. 
\end{pf}

\begin{example} \label{briex}
Let 
$$
\PPP = \{0\text{-dimensional sheaves on } Y \},
$$
and let
$$
\QQQ = \{ E \in \coh (Y) : \Hom (P, E) =0 \text{ for all } P \in \PPP \}.
$$
Lemma~\ref{toda} easily proves that the pair $(\PPP, \QQQ)$ is a torsion pair. 
\end{example}
Let 
\[
\CCC = \coh_{\leq 1}(Y)
\]
denote the full subcategory of coherent sheaves on $Y$ whose support
is of dimension no more than one. We make the following definitions:

$$
\Pp = \{ P \in \CCC | R^{\bullet}\pi_* P \hbox { is a zero-dimensional sheaf on }X \}, 
$$
and
$$
\Qp = \{ F \in \CCC | \hbox{ for all } P \in \Pp, \Hom(P, F) =0 \} = \Pp^{\perp}.
$$
\begin{lemma}\label{lem: (P,Q) is a torsion pair} The pair $(\Pp,
\Qp)$ is a torsion pair in $\CCC$.
\end{lemma}
\begin{pf} 
By lemma~\ref{toda}, it suffices to prove that $\Pp$ is closed under extensions and quotients.

Let $P\pr, P\ppr \in \Pp$, and consider the short exact sequence
$$
0 \to P\pr \to P \to P\ppr \to 0.
$$
We are to show that such a $P$ must live in $\Pp$. Consider now the long exact sequence, 

$$
0 \to \pi_* P\pr \to \pi_* \Pp \to \pi_* P\ppr 
\to R^1\pi_* P\pr \to R^1\pi_* \Pp \to R^1 \pi_* P\ppr \to 0.
$$
Since $P\pr, P\ppr \in \Pp$, we know that $R^1 \pi_*P\pr =0$ and $R^1 \pi_*P\ppr =0$, so $R^1 \pi_*P =0$. We also know that $\pi_*P\pr$ and $\pi_*P\ppr$ are zero-dimensional sheaves, and so it is clear then that $\pi_* P$ must be so as well. This proves that $\Pp$ is closed under extensions. 

Let $P \in \Pp$, and consider a quotient $P\to B \to 0$. Denote the
kernel of this map by $K$. As before, we get a long exact sequence,
$$
0 \to \pi_* K \to \pi_* P \to \pi_* B \to R^1 \pi_* K \to R^1 \pi_* P \to R^1 \pi_* B \to 0.
$$
Since $P \in \Pp, R^1\pi_*P = 0,$ and so $R^1\pi_*B =0$. It remains to
show that $\pi_*B$ is zero-dimensional. We know that $\pi_*K$ is zero
dimensional, since it is a subsheaf of (the zero-dimensional sheaf)
$\pi_*P$. The support of $R^1\pi_*K$ is contained in the singular locus. Suppose $\dim \supp(R^1\pi_*K) =1$. Then, $K$ must have been supported in dimension two, however this contradicts the fact that $K \in \coh_{\leq1}(Y)$.

Hence $R^1 \pi_*K$ is zero dimensional. Now, $\pi_*B$ is  the extension of zero-dimensional sheaves, so it too is zero-dimensional. This completes the proof that $B\in \Pp$, and that $(\Pp, \Qp)$ is a torsion pair.  

\end{pf}

\begin{definition} \label{psp}
A map $\gamma: \OO_Y \to G$ is a $\pi$-\emph{stable pair} if $G\in \Qp$ and $\cok (\gamma) \in \Pp$.\end{definition}
\begin{remark}
Our notion of $\pi$-stable pair is a generalization of the stable pairs of Pandharipande and Thomas \cite{pt1}. In the trivial case when $X=Y$ and $\pi=$ the identity, we have that $(\Pp, \Qp) = (\PPP, \QQQ)$ of example~\ref{briex}, and the \ps$ $ pairs are exactly PT stable pairs.
\end{remark}
\begin{definition} Two $\pi$-stable pairs $\gamma_1:\OO_Y \to G_1$ and $\gamma_2: \OO_Y \to G_2$ are \emph{isomorphic} if there exists a isomorphism of sheaves $\theta: G_1 \to G_2$ making the following diagram commute:
$$
\xymatrix{
\OO_Y \ar[r]^{\gamma_1} \ar[dr]_{\gamma_2} & G_1 \ar[d]^{\theta} \\
 & G_2
}
$$
A \emph{family} of $\pi $-stable pairs on $Y$ over a scheme $T$ is a
coherent sheaf $G$ on $Y\times T$, flat over $T$ and a morphism
$\gamma :\OO _{Y\times T}\to G$ such that for all closed points $t\in
T$, the restriction $\gamma _{t}:\OO _{Y}\to G_{t}$ is a $\pi $-stable
pair.
\end{definition}

\begin{remark} In \cite{bri2}, the tilt of $\AAAA$ with respect to the torsion pair $(\PPP, \QQQ)$ of example~\ref{briex} is denoted $\AAAA^{\#}$, and lemma~2.3 of \cite{bri2} proves that $\AAAA^{\#}$-epimorphisms of the form $\OO_Y \to F$ are precisely stable pairs. The abelian category generated by $\OO_Y$ and $\CCC$ has a tilt whose epimorphisms of the form $\OO_Y \to G$ are precisely $\pi$-stable pairs. This is analogous to the tilt used in \cite{bri2}. However, since it is not strictly necessary for any of our arguments, we will not present a proof here.
\end{remark}

We associate to every $\pi$-stable pair $\gamma:\OO_Y \to G$ a short exact sequence
$$
0 \to \OO_C \to G \to P \to 0,
$$
where $P = \coker (\gamma) \in \Pp$ and $\OO_C = \OO_Y / \ker (\gamma)$.

\begin{prop} \label{prop1}
Let $G$ be a non-zero sheaf. If $\OO_Y \to G$ is a $\pi$-stable pair, then $G$ is not supported exclusively on exceptional curves, and $R^1\pi_*G=0$.
\end{prop}
\begin{pf} Let 
$$
0\to \OO_C \to G \to P \to 0
$$ 
be the associated short exact sequence.  We will first show that $R^1\pi_*\OO_C=0.$ Consider the short exact sequence
$$
0\to I_C \to \OO_Y \to \OO_C \to 0.
$$
Pushing forward yields
$$
0\to \pi_*I_C \to \pi_*\OO_Y \to \pi_*\OO_C \to  \newline
R^1\pi_*I_C \to R^1\pi_*\OO_Y \to R^1\pi_*\OO_C \to 0
$$
which is exact since the dimension of the fibres of $Y\to X$ is at most one. Thus the vanishing of $R^1\pi_* \OO_Y$ by equation~\ref{rationalsing} implies that of $R^1\pi_*\OO_C$. 

Now consider the following long exact sequence.
$$
0\to \pi_* \OO_C \to \pi_* G \to \pi_*P 
\to R^1\pi_* \OO_C \to R^1\pi_* G \to R^1\pi_*P \to 0.
$$ 
From above, we know $R^1\pi_*\OO_C=0$. As well, $P\in \Pp$ implies $R^1\pi_*P=0$. Thus $R^1\pi_*G=0$. 

Now, if $C$ consists of only exceptional curves, then $\pi_*\OO_C$ is zero-dimensional. This implies that $\pi_* G$ is the extension of zero-dimensional sheaves, and therefore zero-dimensional. This means that $G\in \Pp$. By definition of $\pi$-stable pair, $G\in \Qp$. By definition of $\Qp$ the only map from an object of $\Pp$ to an object of $\Qp$ is the zero map, hence the identity of $G$ is the zero map, and $G$ is the zero object.
\end{pf}

Let us introduce some terminology and results taken from \cite{bri2} (modified for our purposes, since we are only interested in sheaves supported in dimension no more than one). Let $\MM$ denote the stack of objects of $\coh_{\leq 1}(Y)$. It is an algebraic stack, locally of finite type over $\CC$. Let $\MM (\OO)$ denote the stack of framed sheaves, that is, the stack whose objects over a scheme $S$ are pairs $(E, \gamma)$ where $E$ is a $S$-flat coherent sheaf on $S\times Y$, of relative dimension no more than one, together with a map $\OO_{S \times Y} \stackrel{\gamma}{\to} E$. Given a morphism of schemes $f: T \to S$, and an object $(F, \delta)$ over $T$, a morphism in $\MM(\OO)$ lying over $f$ is an isomorphism
$$
\theta: f^*(E) \to F
$$
such that the following diagram commutes
 \begin{equation}
 \xymatrix@C=1.5em{ 
 f^*(\OO_{S\times Y}) \ar[d]^{\text{can}} \ar[r]^-{f^*(\gamma)} & f^*(E) \ar[d]^{\theta} \\
\OO_{T\times Y} \ar[r]^{\delta} & F.
}
\end{equation}
The symbol ``$\text{can}$" denotes the canonical isomorphism of pullbacks.

There is a natural map 
\begin{equation} \label{mom}
\MM (\OO) \stackrel{q}{\to} \MM
\end{equation}
sending a sheaf with a section to the underlying sheaf. 

The following lemmas are 2.4 and 2.5 of \cite{bri2}. 
\begin{lemma} \label{bri2.4}
The stack $\MM (\OO)$ is algebraic and the morphism $q$ is representable and of finite type.
\end{lemma}

\begin{lemma} \label{bri2.5} There is a stratification of $\MM$ by locally-closed substacks
$$
\MM_r \subset \MM
$$
such that objects $F$ of $\MM_r(\CC)$ are coherent sheaves satisfying 
$$
\dim_{\CC} H^0(Y, F)=r.
$$
Furthermore, the pullback of the morphism $q$ to $\MM_r$ is a Zariski fibration with fibre $\CC^r$.
\end{lemma}
Both of these are proven in \cite{bri2}.

Let 
\[
\Hilbpbn \subset \MM (\OO )
\]
denote the subcategory of $\MM (\OO )$ consisting of families of $\pi $-stable pairs on $Y$ whose sheaf $G$ has chern character $(0,0,\beta ,n)$.

\begin{lemma} \label{ftalgspace}
$\Hilbpbn $ is a constructible set, that is, it has a finite decomposition into subcategories which are each represented by schemes.
\end{lemma}
\begin{remark}
We expect that $\Hilbpbn $ is in fact represented by a projective
scheme, but we do not pursue that in this paper. The above lemma
suffices for our purposes: the use of $\Hilbpbn $ in the Hall algebra.
\end{remark}
\begin{pf}
Since $\MM (\OO )$ is a locally constructible stack, the content of
the lemma is (1) the subcategory $\Hilbpbn $ is bounded (see
definition~\ref{defn: bounded family}) and (2) the automorphism group
of an object in $\Hilbpbn $ is trivial. We will prove (1) in
Lemma~\ref{bound} and we prove (2) below.

Let $\OO_Y \to G$ be a \ps \hbox{ }pair. We will show that it has only
the trivial automorphism. Consider the associated short exact
sequence,
$$
0 \to \OO_C \to G \to P \to 0.
$$
An automorphism of this \ps \hbox{ }pair leads to a diagram of the
form,
$$
\xymatrix@C=1.5em{ 
0 \ar[r] & \OO_C \ar[d]^{\text{id}} \ar[r]^\gamma & G \ar[d]^g \ar[r]^h & P \ar[r] \ar[d]^{\overline{g}} & 0 \\
0 \ar[r] & \OO_C \ar[r]^\gamma & G \ar[r]^h & P \ar[r] & 0.
}
$$
We will show that $g$ is the identity map. Consider the following
diagram, obtained by subtracting the identity from the diagram above,
\begin{equation} \label{nuffautos}
\xymatrix@C=1.5em{ 
0 \ar[r] & \OO_C \ar[d]_{\text{zero}} \ar[r]^\gamma & G \ar[d]_{g-{\text{ id}}} \ar[r]^h & P \ar@{.>}[dl]_{\delta } \ar[r] \ar[d]^{\overline{g}-{\text id}} & 0 \\
0 \ar[r] & \OO_C \ar[r] & G \ar[r]^h & P \ar[r] & 0.
}
\end{equation}

Since the left-most vertical arrow is zero, a diagram chase proves
that the dotted morphism $\delta $ exists and commutes with the
diagram. However, $P\in \Pp$ and $G \in \Qp$, so $\delta $ must be
zero, and consequently, $\overline{g}-\text{id}=0$. Standard
homological algebra then implies that the morphism $g-\text{id}$ must
be of the form $\gamma \circ \epsilon \circ h$ for some $\epsilon \in
\Hom (P,\OO _{C})$. However, any non-zero $\epsilon $ would give rise
to a non-zero map $\gamma \circ \epsilon :P\to G$ which contradicts
$P\in \Pp$ and $G\in \Qp $. Thus $g-\text{id}=0$ and so $g=\text{id}$.

\end{pf}

\subsection*{The Behrend function identity}
We state and prove a variation of \cite[theorem 3.1]{bri2} of Bridgeland.

\begin{lemma} \label{bflemma} Let $\gamma: \OO_Y \to G$ be a $\pi$-stable pair. Then there is an equality of Behrend's microlocal functions
$$ 
\nu_{\MM (\OO)}(\gamma) = (-1)^{\chi(G)}\nu_{\MM}(G).
$$
\end{lemma}

\begin{pf} The case when $G$ is a stable pair is taken care of by theorem 3.1 of \cite{bri2}. Thus, we may assume that the cokernel $P$ of $\gamma: \OO_Y \to G$ has one-dimensional support. 

Let $\OO_C \sub G$ be the image of $\gamma$. It is the structure sheaf of a subscheme $C\sub Y$ of dimension 1. 
There is a line bundle $L$ on $Y$ such that
\begin{equation} \label{vanish}
H^i(Y, G \otimes L)=0
\end{equation}
for all $i>0$, and there is a divisor $H \in |L|$ such that $H$ meets $C$ at finitely many points, none of which are in the support of $\coker (\gamma)$. 
This claim is verified in lemma~\ref{Jim}. 
 From here, the proof is identical to Bridgeland's, but we include a portion of it to illustrate his ideas.

There is a short exact sequence
$$
0\to \OO_Y \stackrel{s}{\to} L  \to \OO_H(H) \to 0
$$
where $s$ is the section of $L$ corresponding to the divisor $H$. Tensoring it with $G$, and using the above assumptions yields a diagram of sheaves
\begin{equation} \label{diagram}
\xymatrix{
& \OO_Y \ar[d]^{\gamma} \ar[dr]^{\delta} &&& \\
0 \ar[r] & G \ar[r]^\alpha & F \ar[r]^\beta & K \ar[r] &0
}
\end{equation}
where $F= G\otimes L$. The support of the sheaf $K$ is zero-dimensional, and disjoint from the support of $\coker(\gamma)$. In particular, 
\begin{equation} \label{kf}
\Hom_Y(K, F)=0.
\end{equation}

Consider two points of the stack $\MM (\OO)$ corresponding to the maps
$$
\gamma: \OO_Y \to G \text{ and } \delta:\OO_Y \to F.
$$
The statement of lemma~\ref{bflemma} holds for the map $\delta$ because lemma~\ref{bri2.5} together with equation~\ref{vanish} implies that 
$$
q: \MM (\OO) \to \MM
$$
is smooth of relative dimension $\chi (F) = H^0(Y, F)$ over an open neighbourhood of the point $F\in \MM (\CC)$. On the other hand, tensoring sheaves with $L$ defines an automorphism of $\MM$, so the microlocal function of $\MM$ at the points corresponding to $G$ and $F$ are equal. To prove the lemma, it suffices to show that 
$$
(-1)^{\chi (G)}\cdot \nu_{\MM (\OO)}(\gamma) = (-1)^{\chi (F)} \cdot \nu_{\MM(\OO)}(\delta).
$$

Consider the stack $W$ whose $S$-valued points are diagrams of $S$-flat sheaves on $S\times Y$ of the form

$$
\xymatrix{
& \OO_{S \times Y} \ar[d]^{\gamma_S} \ar[dr]^{\delta_S} &&& \\
0 \ar[r] & G_S \ar[r]^{\alpha_S} & F_S \ar[r]^{\beta_S} & K_S \ar[r] &0.
}
$$

There are two morphisms
$$
p: W \to \MM (\OO), \hspace{1cm} q: W \to \MM (\OO),
$$
taking such a diagram to the maps $\gamma_S$ and $\delta_S$ respectively. By passing to an open substack of $W$, we may assume that equation~\ref{kf} holds for all $\CC$-valued points of $W$. It follows that $p$ and $q$ induce injective maps on stabilizer groups of $\CC$-valued points, and hence are representable. 

Recall that Behrend's microlocal function satisfies the property that when $f:T \to S$ is a smooth morphism of relative dimension $d$, there is an identity \cite[proposition 1.5]{beh}
$$
\nu_T = (-1)^df^*(\nu_S).
$$ 
Using this identity, it will be enough to show that at the point $w\in W(\CC)$ corresponding to the diagram~\ref{diagram}, the morphisms $p$ and $q$ are smooth of relative dimension $\chi (K)$ and 0, respectively. For the proof of these facts, see \cite[pages 11--13]{bri2}.
\end{pf}

\begin{lemma} \label{Jim}
Given a $\pi$-stable pair $\gamma: \OO_Y \to G$, we may choose a very ample divisor $H$ on $X$ such that its pull-back is equal to its proper transform (we denote both by $\wt{H}$), it satisfies $\supp (\cok \gamma) \cap \wt{H} =\emptyset$, $\wt{H} \cap \supp G$ is 0-dimensional, and
$$
H^1(Y, G(\wt{H})) =0.
$$
\end{lemma}
\begin{pf} 
First we collect a little notation. Let $E$ be the exceptional locus of $Y$, let $E\pr$ be the image of $E$ in $X$. Define a subset $Z$ as follows
$$
Z = \{ p\in X : G|_{\pi^{-1}(p)} \text{ is one dimensional} \} \subset E\pr.
$$
Notice that $Z$ is a finite collection of points, namely it is the
image under $\pi $ of the exceptional components in the support of
$G$.

Since the cokernel of $\gamma $ is supported in dimension one and it lies in $\Pp $, it is supported on points and exceptional curves. Hence $\pi (\supp (\coker \gamma ))$ is zero dimensional. Moreover, $\pi(\supp G)$ is one
dimensional. Thus we may choose an ample divisor $H$ on $X$ so that
$$
H \cap \pi(\supp(\coker \gamma)) = \emptyset
$$
and $H \cap \pi(\supp G)$ is zero dimensional and does not contain any of the points in $Z$. It follows that $\wt{H} \cap \supp(\coker \gamma)$ is empty, and $\wt{H} \cap \supp (G)$ is zero dimensional. Moreover, by Serre vanishing, we may assume that $H$ is sufficiently ample on $X$ so that 
\begin{equation} \label{Serre}
H^1(X, (\pi_*G)(H))=0.
\end{equation}

We now show that $H^1(Y, G(\wt{H}))=0.$ By proposition~\ref{prop1}, we know that $R^1\pi_*G =0$, since $\OO_Y \to G$ is a $\pi$-stable pair. The sequence
$$
0 \to G \to G(\wt{H}) \to G(\wt{H})|_{\wt{H}} \to 0
$$
gives
$$
\ldots \to R^1\pi_*G \to R^1\pi_*G(\wt{H}) \to R^1\pi_*G(\wt{H})|_{\wt{H}} \to 0.
$$
However, we know $R^1\pi_*G=0$ and $G(\wt{H})|_{\wt{H}}$ is supported on points so $R^1\pi_*G(\wt{H})|_{\wt{H}}=0$, so $R^1\pi_*G(\wt{H})=0$. Now by the Leray spectral sequence, 
$$
\begin{array}{rl}
H^1(Y, G(\wt{H})) & = H^1(X, \pi_*(G(\wt{H}))\\
 & = H^1 (X, \pi_* (G \otimes \pi^*\OO_X(H)))\\
 & = H^1(X, (\pi_*(G))(H)) \\
& =0,
\end{array}
$$
where the last equality comes from equation~\ref{Serre}.
\end{pf}

\section{Stability conditions}

\noindent In this section, we define a stability condition on $\CCC = \coh_{\leq1}Y$. We follow Joyce's treatment of stability conditions as found in section 4 of \cite{joy3}, though not in as great generality.

Let $N_1(Y)$ denote the abelian group of cycles of dimension one modulo numerical equivalence. We begin by quoting lemmas 2.1 and 2.2 of \cite{bri2}. 

\begin{lemma}An element $\beta \in N_1(Y)$ has only finitely many decompositions of the form $\beta = \beta_1 + \beta_2$ with $\beta_i$ effective.
\end{lemma}

\begin{lemma} The Chern character map induces an isomorphism
$$
ch=(\text{ch}_2, \text{ch}_3): F_1K(Y) \to N_1(Y)\oplus \ZZ.
$$
\end{lemma}

Define
$$
\Delta = \{ [E] \in F_1K(Y) : E \in \CCC \}
$$
to be the \emph{positive \emph{or} effective cone} of $F_1K(Y).$

\begin{definition} A \emph{stability condition} on $\CCC$ is a triple $(T, \tau, \leq)$ where $(T, \leq)$ is a set $T$ with a total ordering $\leq$, and $\tau$ is a map $\Delta \stackrel{\tau}{\to} T$ from the effective cone to $T$, satisfying the following condition: whenever $\alpha +\beta = \gamma$ in $\Delta$, then  
$$
\tau (\alpha) < \tau (\gamma) < \tau(\beta),
$$
or
$$
\tau (\beta) < \tau (\gamma) < \tau(\alpha),
$$
or 
$$
\tau (\alpha) = \tau (\gamma) = \tau(\beta).
$$

A triple $(T, \tau, \leq)$ is called a \emph{weak stability condition} if it satisfies the weaker condition that whenever $\alpha +\beta = \gamma$ in $\Delta$, $\tau (\alpha) \leq \tau (\gamma) \leq \tau(\beta)$ or $\tau (\beta) \leq \tau (\gamma) \leq \tau(\alpha)$.
\end{definition}

\begin{definition} A non-zero sheaf $G$ is 
\begin{enumerate}
\item $\tau$-\emph{semistable} if for all $S\subset G$, such that $S \not\iso 0$, we have that $\tau(S)\leq\tau(G/S)$;
\item $\tau$-\emph{stable} if for all $S\subset G$, such that $S \not\iso 0$, we have that $\tau(S) < \tau(G/S)$;
\item $\tau$-unstable if it is not $\tau$-semistable.
\end{enumerate}
\end{definition}

\begin{lemma} \label{basicss} 
Let $F$ and $G$ be $\tau$-semistable sheaves, and let $F \stackrel{f}{\to} G$ be a map of sheaves. Then either $\tau(F) \leq \tau(G)$ or $f=0$.
\end{lemma}
\begin{pf}
Consider the inclusion map $\iota : \img(f) \to G$. Since $G$ is semistable, we know either $\img(f)=0$ or $\tau(\img(f)) \leq \tau(G)$. Consider now the corestriction of $f$, $\cor (f): F \to \img(f)$. Since $F$ is semistable, we know that either $\img (f)=0$ or $\tau(F) \leq \tau (\img (f))$. This implies either $\img(f)=0$ or $\tau (F) \leq \tau (G)$.
\end{pf}

Now we define a stability condition on $\CCC$. 

\begin{definition}\label{defn: mupi stability} Choose an ample divisor $H$ on $X$, let
$\widetilde{H}$ denote the total transform of $H$ in $Y$. Let $A$ be
an ample line bundle on $Y$ and we let $L= \wt{H}+ A$. Note that
$L$ is ample and that $L\cdot C>\wt{H}\cdot C$ for any curve class
$C$.
Given a sheaf $G$ in $\CCC$, define the $\pi$-\emph{slope} of $G$ to be
$$
\mu_{\pi}(G) = \left( \frac{\chi(G)}{\beta \cdot \widetilde{H}}, \frac{\chi(G)}{\beta \cdot L} \right) \in (-\infty, +\infty]\times(-\infty, +\infty],
$$
where by convention $\chi /0=+\infty $ for any $\chi \in \ZZ $, $(-\infty, +\infty]\times(-\infty, +\infty]$ is ordered lexicographically, and $\beta = \beta_G$ the homology class associated to the support of $G$. 
\end{definition}

To make Joyce's Hall algebra machinery work, he introduces the
additional notion of \emph{permissibility} for a (weak) stability
condition \cite[Def.~4.7]{joy3}. The main result of this section is
the following:

\begin{theorem} \label{wpsc}
The map 
$$
\mup : \Delta \to (-\infty, +\infty]\times(-\infty, +\infty]
$$
defines a weak permissible stability condition. 
\end{theorem}

Following definitions 4.1 and 4.7 of \cite{joy3}, we see that we must
prove the following three properties:

\begin{enumerate}
\item (weak seesaw property) for any short exact sequence 
$$
0\to A \to G \to B\to 0,
$$ 
either $\mup(A) \leq \mup(G) \leq \mup(B)$ or $\mup(A) \geq \mup(G) \geq \mup(B),$
\item $\CCC$ is $\mu _{\pi }$-artinian, i.e. there exists no infinite descending chain $\cdots A_2 \subset A_1 \subset A$ in $\CCC$ such that $A_i \neq A_{i+1}$, and $\mu _{\pi }(A_{i+1}) \geq \mu _{\pi }(A_i/A_{i+1})$ for all $i$; and
\item the substack of $\mu _{\pi }$-semistable objects of a fixed
Chern character of the stack parametrizing objects of $\CCC$ is a
constructible substack of $\MM$.
\end{enumerate}
We will prove the first property in lemma~\ref{stab1} and the second
in lemma~\ref{art}. The third property amounts to showing that the
family of $\mu _{\pi }$-semistable sheaves of a fixed chern class is
bounded (see the proof of theorem~4.20 in \cite{joy3}) which we prove
in lemma~\ref{lem: mu-pi semistable sheaves are bounded}.

\begin{lemma} \label{stab1}
The function $\mup$ satisfies the weak seesaw property. 
\end{lemma}
\begin{pf}
Let $0\to A\to G\to B\to 0$ be a short exact sequence of sheaves, and suppose $\mup(A) \leq \mup(G)$. Less concisely, we are supposing
$$
\left( \frac{\chi(A)}{\beta_A \cdot \widetilde{H}}, \frac{\chi(A)}{\beta_A \cdot L}  \right) \leq
\left( \frac{\chi(G)}{\beta_G \cdot \widetilde{H}}, \frac{\chi(G)}{\beta_G \cdot L} \right),
$$
from which we are to deduce that $\mup(G) \leq \mup(B)$. Before we start a case-by-case analysis, notice that $\chi(G) = \chi(A) + \chi(B)$ and $\beta_G = \beta_A + \beta_B$.
\begin{itemize}
\item[case 1:] $\frac{\chi(A)}{\beta_A \cdot \widetilde{H}} < \frac{\chi(G)}{\beta_G \cdot \widetilde{H}}$ and no denominator is zero.

Then this follows from the observation 
$$
\frac{a}{b} < \frac{a+c}{b+d} \Rightarrow \frac{a+c}{b+d} < \frac{c}{d},
$$
provided $b,d>0$. In particular, we assume
$$
\frac{\chi(A)}{\beta_A \cdot \widetilde{H}}\leq
\frac{\chi(G)}{\beta_G \cdot \widetilde{H}}.
$$
Rewriting the second term yields
$$
\frac{\chi(A)}{\beta_A \cdot \widetilde{H}}\leq
\frac{\chi(A)+\chi(B)}{(\beta_A +\beta_B) \cdot \widetilde{H}}.
$$
The observation above then proves that 
$$
\frac{\chi(A)+\chi(B)}{(\beta_A +\beta_B) \cdot \widetilde{H}} \leq
\frac{\chi(B)}{\beta_B \cdot \widetilde{H}},
$$
as desired.
\item[case 2:]$\frac{\chi(A)}{\beta_A \cdot \widetilde{H}} = \frac{\chi(G)}{\beta_G \cdot \widetilde{H}}$, $\frac{\chi(A)}{\beta_A \cdot L} \leq \frac{\chi(G)}{\beta_G \cdot L}$ and no denominator is zero.

We are given that $\frac{\chi(A)}{\beta_A \cdot \widetilde{H}} = \frac{\chi(G)}{\beta_G \cdot \widetilde{H}}$, so
$$
\chi(A) (\beta_G \cdot \widetilde{H}) = (\beta_A\cdot\widetilde{H})\chi(G). 
$$
Writing everything in terms of $A$ and $B$, 
$$
\chi(A) ((\beta_A +\beta_B) \cdot \widetilde{H}) = \beta_A\cdot\widetilde{H}(\chi(A)+\chi(B)),
$$
which implies
$$
\chi(A)(\beta_B\cdot \widetilde{H}) = (\beta_A\cdot \widetilde{H})\chi (B).
$$
Since we assume that all denominators are non-zero, we have
$$
\frac{\chi(A)}{\beta_A\cdot \widetilde{H}} = \frac{\chi(B)}{\beta_B\cdot \widetilde{H}}.
$$
So we must show that 
$$
\frac{\chi(G)}{\beta_G \cdot L} < \frac{\chi(B)}{\beta_B \cdot L}.
$$
This follows from the same observation made in case 1.

\item[case 3:] $\beta_A \cdot \widetilde{H} = 0$.  Then $+\infty =
\frac{\chi(A)}{\beta_A\cdot\widetilde{H}} \leq \frac{\chi(A)
+\chi(B)}{(\beta_A + \beta_B)\cdot \widetilde{H}}$. This implies that
$\frac{\chi(A) +\chi(B)}{(\beta_A + \beta_B)\cdot \widetilde{H}} =
+\infty$, so $(\beta_A + \beta_B)\cdot \widetilde{H}=0$, and hence
$\beta_B \cdot \widetilde{H}=0$. This reduces us to Gieseker stability
on $Y$, which we know satisfies the weak seesaw property.

\item[case 4:] $\beta_G\cdot \widetilde{H} = 0$.  We know $\beta \cdot
\pi^*H = \pi_*(\beta\cdot\pi^*H) = \pi_*(\beta) \cdot H \geq 0$, hence
$\beta \cdot \widetilde{H} \geq 0$ for any effective curve class
$\beta$, so we must have $\beta_A \cdot\widetilde{H} = 0$ and $\beta_B
\cdot\widetilde{H} = 0$. This lands us back in the case of Gieseker
stability on $Y$, and lemma~\ref{stab1} is proven.
\end{itemize}
\end{pf}

\begin{lemma} \label{art}
The category $\CCC$ is $\mup$-artinian.
\end{lemma}
\begin{pf}
Joyce proves that a weak stability condition is Artinian if it is
dominated by an Artinian weak stability condition
\cite[4.10,4.11]{joy3}. Recall that a weak stability condition
$\tilde{\tau}$ is said to dominate $\tau$ if for any $A,B$ in $\CCC$
with $\tau(A) \leq \tau(B)$ then $\tilde{\tau}(A) \leq
\tilde{\tau}(B)$.  Let
$$
\delta(G) = -\dim \supp G \in \ZZ,
$$
then $\delta $ is an Artinian, weak stability condition \cite[4.19]{joy3}.
Thus to prove the lemma, it suffices to show that $\mup$ is dominated by $\delta$.

Let $\mup(A) \leq \mup(B)$. We need to show that this implies that $\delta(A) \leq \delta(B)$. Expanding, we have that
$$
\left(\frac{\chi(A)}{\beta_A \cdot \widetilde{H}}, \frac{\chi(A)}{\beta_A \cdot L}\right) \leq \left( \frac{\chi(B)}{\beta_B \cdot \widetilde{H}}, \frac{\chi(B)}{\beta_B \cdot L} \right).
$$
We proceed with a case-by-case analysis. 
\begin{itemize}
\item[case 1:] the denominators are non-zero. Since we have restricted our attention to sheaves supported in dimension $\leq$ 1, it follows that neither $A$ nor $B$ is $0$-dimensional. Thus, they are both one dimensional, and $\delta (A) = \delta (B)$. In particular, $\delta (A) \leq \delta (B)$. 

\item[case 2:] $\beta_B \cdot \widetilde{H} = 0 $ and $\beta_A \cdot \widetilde{H} \neq 0$. Then $\dim \supp A \geq 1 \geq \dim \supp B$. So $\delta (A) \leq \delta (B)$.

\item[case 3:] Both $\beta_B \cdot \widetilde{H} = 0 $ and $\beta_A \cdot \widetilde{H} = 0$. Then $\mup (A) \leq \mup (B)$ amounts to regular Geiseker stability, which is dominated by $\delta$ as demonstrated by Joyce \cite[\S 4.4]{joy3}.

\end{itemize}
\end{pf}

To finish the proof that $\mup$ is a permissible weak stability
condition, it remains only to prove the family of all
$\mup$-semistable sheaves of a fixed Chern class is bounded. This is
proven in lemma~\ref{lem: mu-pi semistable sheaves are bounded}  in the next section.

\subsection*{Boundedness}

In this section, we prove that the family of $\pi $-stable pairs with
fixed chern classes is bounded (lemma~\ref{bound}) and we prove that
the family of $\mu _{\pi }$-semistable sheaves with fixed chern
classes is bounded (lemma~\ref{lem: mu-pi semistable sheaves are
bounded}).

We begin by recalling some basic results concerning boundedness (cf. \cite{hule}). 

\begin{definition}\label{defn: bounded family} A subcategory $\UU$ of
$\coh (Y)$ is bounded if there exists a scheme $S$ of finite type and
a sheaf $U$ on $X\times S$ such that for every object $U_i$ of $\UU$,
there exists a closed point $s_i \in S$ such that $U_i \iso U|_{X
\times \{s_i\}}$.
\end{definition}
Notice that this definition still makes sense if we have a set of isomorphism classes of sheaves instead of a category.

\begin{definition} Let $Y$ be a scheme, let $\OO(1)$ be an ample line bundle, and let $m$ be an integer. A sheaf $F$ on $Y$ is $m$-regular if, for all $i>0$, 
$$ 
H^i(Y, F(m-i))=0.
$$
\end{definition}
\indent A proof for the following may be found in \cite{kle}, as well as in \cite{mum}.
\begin{lemma} If $F$ is $m$-regular, then the following statements are true:
\begin{enumerate} 
\item $F$ is $m\pr$-regular for all $m\pr \geq m$.
\item $F(m)$ is globally generated.
\item For all $n\geq 0$, the natural map $H^0(Y, F(m)) \otimes H^0(Y, \OO(n)) \to H^0(Y, F(n+m))$ is surjective.
\end{enumerate}
\end{lemma} 

\begin{definition} The Mumford-Castelnuovo regularity of a sheaf $F$ is the number $\reg(F) = \inf \{m \in \ZZ : F \mbox{ is } m\mbox{-regular } \}.$
\end{definition}

\begin{lemma} Let $\UU$ be a category of sheaves on $Y$. The following statements are equivalent. 
\begin{enumerate}
\item $\UU$ is bounded.
\item The set of Hilbert polynomials of objects $U_i$ of $\UU$ is finite, and there is an integer $N$ such that for all objects $U_i$ of $\UU$, $\reg U_i < N$.
\item The set of Hilbert polynomials of objects $U_i$ of $\UU$ is finite, and there exists a sheaf $F$ such that each object of $\UU$ is isomorphic to a quotient of $F$.
\end{enumerate}
\end{lemma}
The proof of this lemma may be found in \cite{gro}.

\begin{lemma} \label{bound}
 The family of $\pi$-stable pairs $\OO_Y \stackrel{\gamma}{\to} G$
with a fixed Chern class is bounded.
\end{lemma}

\begin{pf} To each such $\pi$-stable pair there is an associated a short exact sequence, 
$$
0\to \OO_C \to G \to P \to 0,
$$
where $\OO_C$ is the image of the map $\gamma$, and $P$ is the
cokernel. We will show that the family of possibilities for $\OO_C$
and the family of possible $P$s are both bounded families. Once this
is established, it is clear that the family of sheaves underlying a
$\pi$-stable pair is a bounded family of sheaves.

First we will consider the family of possibilities for $\OO_C$. To show that this family is bounded, we will show:
\begin{enumerate}
\item the Hilbert polynomials of this family take only a finite number of values; and
\item there exists a single sheaf that surjects onto each member of this family. 
\end{enumerate}
The second requirement is trivially satisfied, since each member of
this family is the structure sheaf of a subscheme of $Y$, and hence,
admits a surjective map from $\OO_Y$. It remains to find upper and
lower bounds for the coefficients of the Hilbert polynomial of a
general element from this family.

In contrast to $PT$-theory, the support of $\OO_C$ is not equal to the support of $G$, since $P$ is not necessarily zero-dimensional. However, we still have $\beta_G = \beta_C + \beta_P$, where all $\beta$ are effective. We know that there are only finitely many decompositions of $\beta_G$ into the sum of two effective curve classes. This forces an upper and lower bound on the linear coefficient of the Hilbert polynomial of $\OO_C$. It remains to find upper and lower bounds for the Euler characteristic of $\OO_C$ (the constant coefficient of the Hilbert polynomial). 

The Leray spectral sequence proves that $\chi(P) = \chi(R^\bullet\pi_*P)\geq0$, the inequality following from the fact that $R^\bullet\pi_*P$ is zero dimensional. Now, $\chi(G) = \chi (P) + \chi (\OO_C)$, and $\chi(P)\geq0$ implies $\chi (\OO_C) \leq \chi (G)$. This gives us an upper bound, since the Chern character, and hence, the Euler characteristic, of $G$ is fixed. For the lower bound, let $\alpha_G = (\beta_G, n_G)$ be the Chern character of $G$, and let $\alpha_C = (\beta_C,n_C)$ be the Chern character of $C$. 

In general, if $\Hilb^{(\beta, n)}$ is non-empty (say one of its points represents a curve $J$), then $\dim \Hilb^{(\beta, n+k)} \geq 3k$ since we get a $3k$-dimensional space of curves coming from the curve $J$ with $k$ ``wandering points." This line of reasoning tells us that 
$$
\dim \Hilb^{(\beta_C,n_G)} \geq 3(n_G -n_C).
$$
Rearranging this yields
$$
n_C \geq n_G - \frac{1}{3}\dim \Hilb^{(\beta_C, n_G)} .
$$
This gives us a lower bound for $n _C$, which completes the proof that the corresponding family is bounded. 

Now to show that the family of cokernels is bounded, we will show
\begin{enumerate}
\item the Hilbert polynomials of this family take on only a finite number of values; and
\item there is a common upper-bound to the index of regularity.
\end{enumerate}

Using the Leray spectral sequence again, we note that for all $P\in\Pp$,
$$
H^1(Y, P) = H^1(X, \pi_*P)=0,
$$
thus, all $P\in\Pp$ are 1-regular. To show this family is bounded, it remains to find upper and lower bounds for the coefficients of the Hilbert polynomial of a general object.

As above, there are only a finite number of options for the support curve of $P$. This yields upper and lower bounds on the linear coefficient of the Hilbert polynomial. 

We know that $\chi (G) = \chi (P) + \chi (\OO_C)$. Since $\chi (\OO_C)$ is bounded, and $\chi (G)$ is fixed, so too must $\chi(P)$ be bounded. 

This completes the proof that the sheaves underlying a $\pi$-stable
pair of fixed $K$-class forms a bounded family of sheaves.
\end{pf}

\begin{lemma}\label{lem: mu-pi semistable sheaves are bounded}
The family of $\mu _{\pi }$-semistable sheaves with fixed chern
character $(0,0,\beta ,n)$ is bounded.
\end{lemma}
\begin{pf}
For a sheaf $G$ of dimension one, we use the notation $\beta _{G} $ to
denote the corresponding curve class and we let
\[
\mu _{N} (G) = \frac{\chi (G)}{N\cdot \beta _{G}}\in (\infty ,\infty ]
\]  
be the $N$-slope, for any $\QQ $-divisor $N$. Note that $\mu _{\pi }$
still denotes the $\pi $-slope so that in this notation
\[
\mu _{\pi } (G)= (\mu _{\wt{H}} (G),\mu _{L} (G))\in (\infty ,\infty ]\times (\infty ,\infty ].
\]
We will construct an ample divisor $A_{\epsilon }$ such that every
$\pi $-semistable sheaf $F$ of chern character $(0,0,\beta ,n)$ is
either $\mu _{L}$-semistable or $\mu _{A_{\epsilon }}$-semistable. The
lemma will then follow since for any ample divisor $N$, the family of
$\mu _{N}$-semistable sheaves of fixed chern classes form a bounded
family \cite[Thm~3.3.7]{hule}.

Let $F$ be a $\mu _{\pi }$-semistable sheaf with $ch (F)= (0,0,\beta
,n)$. We may assume that $ \wt{H}\cdot \beta >0$ since if $\wt{H}\cdot
\beta =0$, then the $\mu _{\pi }$-semistability of $F$ implies $\mu
_{L}$-semistability and we are done.

We construct our ample $A_{\epsilon }$ as follows. Let $A$ be an ample $\QQ $-divisor with $A\cdot \beta =\wt{H}\cdot \beta $ and let 
\[
A_{\epsilon } = (1-\epsilon )\wt{H} +\epsilon A.
\]
Since $\wt{H}=\pi ^{*} (H)$ is in the boundary of the nef cone and $A$
is ample, $A_{\epsilon }$ is ample for any $\epsilon \in \QQ \cap
(0,1)$. We note that $A_{\epsilon }\cdot \beta =\wt{H}\cdot \beta $
for all $\epsilon $. We will choose an appropriate $\epsilon $ below.

Since there are a finite number of decompositions $\beta =\beta
_{1}+\beta _{2}$ with $\beta _{i}$ effective \cite[Lemma~2.1]{bri2}, the set
\[
\{\wt{H}\cdot \beta _{K} \}_{K\subset F}
\]
is finite. By $\mu _{\pi }$-semistability, we know that 
\[
\mu _{\wt{H}} (K)\leq \mu _{\wt{H}} (F)
\]
for all $K\subset F$. Since the set $\{\wt{H}\cdot \beta _{K}
\}_{K\subset F}$ is finite, there exists some $\delta >0$ such that
\[
\mu _{\wt{H}} (K)+\delta < \mu _{\wt{H}} (F)
\]
for all $K\subset F$ such that $\mu _{\wt{H}} (K)<\mu _{\wt{H}}
(F)$. In other words, if the $\wt{H} $-slope of a subsheaf $K\subset
F$ is strictly less than the $\wt{H}$-slope of $F$, then it is bounded
away from the $\wt{H}$-slope of $F$ by $\delta $, a number independent
of $K$ and $F$ (but depending on $\beta $ and $n$). If this were not
the case, there would have to be an infinite number of possible
denominators in $\mu _{\wt{H}} (K)=\chi (K)/ (\wt{H}\cdot \beta _{K})$
which is not true.

We now choose $\epsilon >0$ small enough so that 
\[
\epsilon \cdot \mu _{\wt{H}} (F)\cdot \left(1-\frac{A\cdot \beta _{K}}{\wt{H}\cdot \beta _{K}} \right) <\delta 
\]
for all $K\subset F$ with $\mu _{\wt{H}} (K)<\mu _{\wt{H}} (F)$. Then
for all such $K\subset F$ we get
\[
\epsilon \cdot \mu _{\wt{H}} (F)\cdot \left(1-\frac{A\cdot \beta
_{K}}{\wt{H}\cdot \beta _{K}} \right) +\mu _{\wt{H}} (K) <\delta +\mu
_{\wt{H}} (K)<\mu _{\wt{H}} (F)
\]
which implies
\[
\epsilon\, \frac{\chi (F)}{\wt{H}\cdot \beta }\,\, -\,\, \epsilon\, \frac{\chi (F) (A\cdot \beta _{K})}{(\wt{H}\cdot \beta ) (\wt{H}\cdot \beta _{K})} \,\,+\,\,\frac{\chi (K)}{\wt{H}\cdot \beta _{K}} <\frac{\chi (F)}{\wt{H}\cdot \beta }.
\]
Clearing denominators and rearranging, we get 
\begin{align*}
\chi (K)\, \wt{H}\cdot \beta &\,<\,\chi (F)\,(\epsilon \,A\cdot \beta _{K}+ (1-\epsilon )\wt{H}\cdot \beta _{K} ) \\
&\,=\, \chi (F)\,(A_{\epsilon }\cdot \beta _{K}).
\end{align*}
Using the fact that $A_{\epsilon }\cdot \beta =\wt{H}\cdot \beta $ the above implies
\[
\frac{\chi (K)}{A_{\epsilon }\cdot \beta _{K}} < \frac{\chi (F)}{A_{\epsilon }\cdot \beta }
\]
So we've proved that 
\[
\mu _{A_{\epsilon }} (K)<\mu _{A_{\epsilon }} (F)
\]
for all $K\subset F$ with $\mu _{\wt{H}} (K)<\mu _{\wt{H}} (F)$.

This is now enough to prove our claim: if $F$ is $\mu _{\pi
}$-semistable with $ch (F)= (0,0,\beta ,n)$, then either $F$ is $\mu
_{A_{\epsilon }}$-semistable or $\mu _{L}$-semistable, for if not,
then there exists $K\subset F$ such that $\mu _{A_{\epsilon }} (K)>\mu
_{A_{\epsilon }} (F)$ and $\mu _{L} (K)>\mu _{L} (F)$. But then $\mu
_{\pi } (K)\leq \mu _{\pi } (F)$ implies $\mu _{\wt{H}} (K)<\mu
_{\wt{H}} (F)$ which then by construction implies $\mu _{A_{\epsilon
}} (K)<\mu _{A_{\epsilon }} (F)$ which is a contradiction.

Thus the family of $\mu _{\pi }$-semistable sheaves of chern character
$(0,0,\beta ,n)$ is contained in the union of the families of $\mu
_{A_{\epsilon }}$-semistable and $\mu _{L}$-semistable sheaves of
chern character $(0,0,\beta ,n)$ and is thus bounded.   
\end{pf}

\section{The torsion pair and the stability condition}

In this section, we show that $\Pp$ may be conveniently expressed in terms of the stability condition, and similarly for $\Qp$. First we give a rapid introduction to the modern Harder-Narasimhan property, a generalization of the Harder-Narasimhan filtration of \cite{hana}.

\begin{definition} A weak stability condition $(T, \tau, \leq)$ on $\CCC$ is said to have the \emph{Harder-Narasimhan property} if for every sheaf $G$, there exists a unique filtration of $G$
$$
0 = HN_\tau(G)_0 \subset HN_\tau(G)_1 \subset \cdots \subset HN_\tau(G)_{N-1} \subset HN_\tau(G)_N=G
$$
(where the inclusions are strict) such that the quotients 
$$
Q_i = HN_\tau(G)_i / HN_\tau(G)_{i-1}
$$ 
are $\tau$-semistable and 
$$
\tau (Q_i) > \tau (Q_{i+1})
$$
for all $i>0$. When it is clear from the context, most of the notation will be suppressed, and we will denote the Harder-Narasimhan filtration of $G$ with respect to $\tau$ by $0\subset G_1 \subset G_2 \subset \cdots \subset G_{N-1} \subset G_N=G$. The $G_i$ are called the \emph{filtered objects} of the Harder-Narasimhan filtration, and the $Q_i$ are called the \emph{quotient objects}. 
\end{definition}

We borrow the following definition and theorem from Joyce \cite{joy3}. In \cite[\S 9]{joy1}, Joyce proves that the category of coherent sheaves satisfies assumptions 3.7 of \cite{joy3}. This is enough for us to conclude that the assumptions are also true of $\CCC$ the category of coherent sheaves supported in dimension one or less.

\begin{theorem}[\cite{joy3}, Theorem 4.4] Let $(T, \tau, \leq)$ be a weak stability condition on an abelian category $\AAAA$. If $\AAAA$ is Noetherian and $\tau$-artinian, then $(T, \tau, \leq)$ has the Harder-Narasimhan property.
\end{theorem}

\begin{corollary} The weak stability condition $\mup$ on $\CCC$ has the Harder-Narasimhan property. 
\end{corollary}
\begin{pf} The category $\CCC$ is Noetherian because it is a subcategory of the category of coherent sheaves, which is Noetherian. Corollary~\ref{art} proves that $\CCC$ is $\mup$-artinian.
\end{pf}

When we refer to the Harder-Narasimhan filtration in what follows, we will always be referring to the filtration with respect to the stability condition $\mup$.

We present some notation before we state and prove the main result of this section. Recall that our slope function $\mup$ takes values in the lexicographically ordered set $(-\infty, +\infty]\times (-\infty, +\infty].$ To avoid awkwardly writing the ordered pairs $(+\oo, +\oo)$ and $(+\oo, 0)$ through-out, let us denote 
$$
\oo := (+\infty, +\infty),
$$
and 
$$
\twoo:= (+\infty, 0).
$$
Given an interval $I \sub (-\infty, +\infty]\times (-\infty, +\infty]$, we define $\SSS( I )\sub \CCC$ to be the full subcategory of zero objects together with those one-dimensional sheaves whose Harder-Narasimhan quotients have $\mup$-value in the interval $I$. If $a,b\in (-\infty, +\infty]\times (-\infty, +\infty]$ such that $a<b$, then we denote the closed interval between $a$ and $b$ by $a\leq \square \leq b$, and similarly for open, half-open, etc. intervals.

\begin{lemma} \label{pqss}
$$
\Pp = \SSS(\square \geq \twoo),
$$ 
and
$$
\Qp = \SSS( \square < \twoo ).
$$

\end{lemma}

\begin{pf}
First we will show that $\Pp \subset \SSS (\square \geq \twoo).$

Case 1: let $P \in \Pp$ be semi-stable. We will show that $P \in \SSS (\square \geq \twoo)$.
Since $P$ is semi-stable, it suffices to show that $\mup (P) \geq \twoo$. Now $P \in \Pp$ implies that $\chi (P) \geq 0$. By ampleness of $L$, we know $\beta_P \cdot L \geq 0$. Hence $\mup (P) \geq (+\infty, 0)$.

Case 2: Let $P \in \Pp$ be general, let the following be its Harder-Narasimhan (HN) filtration, 
$$
0 = P_0 \subset P_1 \subset \cdots \subset P_{N-1} \subset P_N = P,$$
and let $Q_i = \frac{P_i}{P_{i-1}}$ be the $i$th quotient; we must show that $\mup (Q_i) \geq \twoo$. Since $P_N = P \in \Pp$ and $\Pp$ is closed under quotients, it follows that $Q_N \in \Pp$. By definition of the HN filtration, $Q_N$ is semi-stable, hence $\mup (Q_N) \geq \twoo$ and $Q_N \in \SSS (\square \geq \twoo)$ by the previous case. Another defining property is that $\mup (Q_1) > \mup (Q_2) > \cdots > \mup (Q_N)$. Hence, $\mup (Q_i) \geq \twoo$ for all $i$, in other words, $P \in \SSS (\square \geq \twoo)$.

Now we will show that $\SSS (\square \geq \twoo) \subset \Pp$. 

Case 1: Let $G \in \SSS (\square \geq \twoo)$ be semi-stable. In this
case, $G \in \SSS (\square \geq \twoo)$ implies that $\mup (G) \geq
\twoo$. Since $(\Pp, \Qp)$ is a torsion pair, for every sheaf there
exists a uniquely associated short exact sequence,
$$
0 \to A \to G \to B \to 0
$$
where $A\in \Pp$ and $B\in\Qp$. By the above, we know that $A \in \SSS
(\square \geq \twoo)$. By the semi-stability of $G$, $\mu _{\pi }
(B)\geq \mu _{\pi } (G)\geq \twoo $ and hence we have $\chi (B) \geq
0$.

Notice that $G$ must be supported on a fibre of $\pi$ because $\mup
(G) \in \{+\infty\}\times (0, +\infty]$. Hence $B$ is also supported
on a fibre. We claim that this forces $H^0(B) =0$. For suppose there
was a non-zero map $\OO_Y \to B$. This would yield non-trivial $0\to
\OO_C \to B$ where $C$ is the support of the map $\OO_Y \to B$.  From
the proof of proposition~\ref{prop1}, we know that $R^1\pi_\ast \OO_C
=0$, which implies that $\OO_C \in \Pp$ which contradicts the
definition of $\Qp$. Hence $H^0(B)=0$.

However, $\chi(B) \geq 0$ so $\dim H^1(B) \leq 0$. This implies that
$H^1 (B) =0$, which implies that $R^1\pi_*B =0$ (by the theorem of
cohomology and base-change) and hence $B\in \Pp$. Since $B\in\Qp$ we
conclude that $B=0$ and $G=A \in \Pp$.

Case 2: Let $G\in \SSS (\square \geq \twoo)$ be general. We need to
show that $G \in \Pp$. Let
$$
0 = G_0 \subset G_1 \subset \cdots \subset G_{N-1} \subset G_N = G
$$
be the HN filtration of $G$, and let $Q_i = \frac{G_i}{G_{i-1}}$
denote the corresponding semistable quotients. By assumption, $\mup
(Q_i) \geq \twoo$; notice that $G_1 = Q_1$, so we have that $G_1$ is
semistable and $\mup (G_1) \geq \twoo$. The previous case then proves
that for all $i$, $Q_i\in \Pp$, and since $\Pp $ is closed under
extensions, we see that $G\in \Pp$.

This completes the proof that $\Pp = \SSS (\square \geq \twoo)$. 


The pairs $(\Pp ,\Qp )$ and $(\SSS (\square \geq \twoo ),\SSS (\square
< \twoo ))$ both form torsion pairs; the former we proved in
Lemma~\ref{lem: (P,Q) is a torsion pair}, the later because $\mu _{\pi
}$ is a stability condition. Since any torsion pair is completely
determined by its torsion part, $\Qp =\SSS (\square <\twoo )$ follows
from $\Pp =\SSS (\square \geq \twoo )$ and the lemma is proved.
  
\end{pf} 

\section{The motivic Hall algebra}
Here we provide a quick summary of the constructions and results of Bridgeland's papers \cite{bri1}, \cite{bri2} (which came into existence as a gentle introduction to part of Joyce's theory of motivic Hall algebras \cite{joy1}, \cite{joy2}, \cite{joy3}, \cite{joy4} ).

Let $S$ be a stack, locally of finite type over $\CC$ and with affine stabilizers. 

\begin{definition} \label{hall} The \emph{relative Grothendieck group $K(\St/S)$ of stacks over $S$} is the $\QQ$-vector space spanned by symbols
$$
\left[ T \stackrel{m}{\to} S \right]
$$
(where $T$ is a finite-type stack and $m$ is a morphism), subject to the following relations. 
\begin{enumerate}
\item[a)] $\left[ T\to S\right] = \left[ U\to S \right] + \left[ F \to S \right]$ where $U$ is an open substack of $T$ and $F$ is the corresponding closed complement.
\item[b)] $\left[T_1 \stackrel{s\circ f}{\to} S \right] = \left[ T_2 \stackrel{s\circ g}{\to} S \right]$, if $T_1 \stackrel{f}{\to} B$ and $T_2 \stackrel{g}{\to} B$ are Zariski fibrations\footnote{a Zariski-local product space} over $B$ with identical fibres, and $B\stackrel{s}{\to} S$ is a morphism of stacks.
\item[c)] $\left[T \stackrel{a}{\to} S \right] = \left[ T\pr \stackrel{b}{\to} S \right]$ if there exists a commutative diagram
\begin{equation*} \xymatrix@C=1em{
T\ar[rr]^c \ar[dr]^a &  & T\pr\ar[dl]_b \\  &S }
\end{equation*}
such that the associated map on $\CC$-points $T(\CC) \stackrel{c}{\to} T\pr(\CC)$ is an equivalence of categories.
\end{enumerate}
\end{definition}

The vector space $K(\St / S)$ is a $K(\St / \spec \CC)$-module, whose action we now describe. Let $\left[ A \to \spec \CC\right] = \left[ A \right] \in K(\St / \spec \CC)$, and let $\left[T \stackrel{m}{\to} S\right] \in K(\St / S)$. Then we define 
$$
\left[A \right] \cdot \left[T \stackrel{m}{\to} S\right] = \left[ A \times_{\spec \CC}T \stackrel{f}{\too} S \right]
,$$
where $f$ is the composition of the projection $A \times_{\spec \CC}T \to T$ and the map $T \to S$.

We are most interested in the case where $S$ is the stack of objects
in $\CCC $, i.e. coherent sheaves on $Y$ supported in dimension one or
less. We denote this stack by $\MM$, and we denote $K (\St / \MM)$ by
$\HA$. The vector space $\HA$ is the motivic Hall algebra; let us
justify the name by endowing it with the structure of an
algebra. First, we define $\MM^{(2)}$ to be the stack of short-exact
sequence of sheaves on $\MM$. Now, given $\left[ A \to \MM \right]$ and
$\left[ B \to \MM \right]$ we define the convolution product $\left[ A \to
\MM \right] * \left[ B \to \MM \right]$ to be $\left[ Z \to \MM\right]$ where
$Z\stackrel{c\circ g}{\to} \MM$ is defined by the following Cartesian
diagram:
\[
\begin{CD}
Z & @>g>> &\MM^{(2)}&@>c>>&\MM \\
 @VVV & &@VV(l,r)V \\
A \times B&@>>> &\MM\times\MM.
\end{CD}
\] 
The morphisms $l,r$ are the ``left hand" and ``right hand" morphisms, which project a short exact sequence to its left-most (resp. right-most) non-zero entry. The morphism $c$ is the ``centre" morphism. Intuitively, given families of sheaves $A \to \MM$ and $B \to \MM$ their product in the Hall algebra is the family $Z \to \MM$ parametrizing extensions of objects of $B$ by objects of $A$.

The motivic Hall algebra is useful tool. It holds enough information to allow us to retrieve Euler characteristics, yet is flexible enough to produce decompositions of elements in terms of extensions. We will describe an ``integration" map on $\HA$ taking values in the ring of polynomials. Equations among elements of $\HA$ will be integrated to yield equations of polynomials. This entire framework will then be souped-up to incorporate Laurent series, and our theorem will be the result of applying the souped-up integration map to equations in the souped-up Hall algebra.

In \cite{bri1}, Bridgeland introduces  \emph{regular elements}. Let $\Kvar$ denote the \emph{relative Grothendieck group of varieties over} $\CC$ (cf. definition~\ref{hall}). Let $\LL$ denote the element $[ \AAA^1 \to \CC]$, the Tate motive. Consider the maps of commutative rings, 
$$
\Kvar \to \Kvar [\LL^{-1}] \to \KSt,
$$
and recall from \cite{bri1} that $\HA$ is an algebra over
$\KSt$. Define a $\Kvar [\LL^{-1}]$-module
$$
\HRA \sub \HA
$$
to be the span of classes of maps $[V \stackrel{f}{\to} \MM]$ with $V$ a variety. We call an element of $\HA$ \emph{regular} if it lies in this submodule. The following result is theorem 5.1 of \cite{bri1}.

\begin{theorem}
The submodule of regular elements is closed under the convolution product:
$$
\HRA \ast \HRA \sub \HRA,
$$
and is therefore a $\Kvar [\LL^{-1}]$-algebra. Moreover the quotient 
$$
\HSC = \HRC / (\LL -1)\HRC,
$$
is a commutative $\Kvar$-algebra.
\end{theorem}
Bridgeland equips $\HSC$ with a Poisson bracket, defined by 
$$
\{ f, g \} = \frac{f*g - g*f}{\LL -1}.
$$
The integration map $I$ is defined on $\HSC$. Now we work toward the polynomial ring in which it takes values.

Recall that $K(Y)$ is the numerical $K$-theory of $Y$. Recall $\Delta \subset F_1K(Y)$ is the effective cone of $F_1K(Y)$, that is, the collection of elements of the form $[F]$ where $F$ is a one-dimensional sheaf. Define a ring $\CC[\Gamma]$ to be the vector space spanned by symbols $x^{\alpha}$ for $\alpha \in \Delta$ and defining the multiplication by 
$$
x^{\alpha} \cdot x^{\beta} = x^{\alpha + \beta}.
$$
We equip $\CC [\Delta]$ with the trivial Poisson bracket. We are now ready for the following theorem. 

\begin{theorem}[5.1 of \cite{bri2}] There exists a Poisson algebra homomorphism
$$
I: \HSC \to \CC[\Delta]
$$
such that 
$$
I( \left[Z \stackrel{f}{\to} \MM_{\alpha}\right] ) = \chi( Z, f^*(\nu) )x^{\alpha}
,$$
where $\nu: \MM \to \ZZ$ is Behrend's microlocal function of $\MM$, and $\MM_\alpha$ denotes the component of $\MM$ with fixed Chern character $\alpha$.
\end{theorem}


\section{Equations in the infinite-type Hall algebra and the fake proof}

For the sake of exposition only, we follow \cite{bri2} and \cite{cale} by introducing an infinite-type version of the Hall algebra. This has the benefit of allowing non-finite-type stacks, but the devastating draw-back of not admitting an integration map. We use it because it will allow us to temporarily work without having to think about convergence of power series. Also, many of the arguments will be used again later. We end this chapter with a fake proof of our main result. It is our hope that this fake proof helps the reader to navigate the true one in the following chapter.

The \emph{infinite-type Hall algebra} is defined by considering symbols as in definition~\ref{hall}, but with $T$ assumed only to be \emph{locally} of finite type over $\CC$, and use relations as before, except that we do not use relation (a). (Admitting relation (a) in this case would make every infinite-type Hall algebra trivial). We denote it by $\HHHH _\infty (\CCC)$. 

Given a substack $\NN \sub \MM$, we let
$$
1_{\NN} = [\NN \stackrel{i}{\to} \MM]
$$
denote the inclusion $i:\NN \to \MM$. Pulling back the morphism (\ref{mom}) to $\NN \sub \MM$ gives a stack denoted $\NN (\OO)$ with a morphism $\NN (\OO) \stackrel{q}{\to} \NN$, and hence an element
$$
1_{\NN}\!\!\!\!\!^{\OO} = [\NN (\OO) \stackrel{q}{\to} \MM] \in \HHHH _\infty (\CCC).
$$

For example, $\Pp$ and $\Qp$ are full subcategories of $\CCC$, and define substacks of $\MM$, which we abusively denote with the same letters, $\Pp, \Qp \sub \MM$. These substacks define elements of the infinite-type Hall algebra, 
$$
\onePp, \oneQp\in \HHHH_\infty (\CCC).
$$
Other examples include 
$$ 
\HHH = [ \Hilb (Y) \to \MM],
$$ 
$$
\HHe = [ \Hilbe \to \MM ],
$$ and 
$$
\HHp = [ \Hilbp(Y) \to \MM] \in \HHHH_\infty (\CCC),$$ 
where $\Hilbe$ denotes the Hilbert scheme of curves supported on fibres of $\pi$, and the map to $\MM$ is given by taking $\OO_Y \to G$ to $G$. Note that all Hilbert schemes are restricted to the components parametrizing sheaves $G$ of dimension one.

\begin{lemma} \label{basiclemma}
$$
1_{\CCC} = \onePp \ast \oneQp
$$
\end{lemma}
This lemma reflects the fact that $(\Pp, \Qp)$ is a torsion pair.

\begin{pf} 
Form the following Cartesian diagram:
\begin{equation} \label{ugh}
\begin{CD}
Z & @>f>> &\MM^{(2)}&@>b>>&\MM \\
 @VVV & &@VV(a_1,a_2)V \\
\Pp\times\Qp &@>i>> &\MM\times\MM.
\end{CD}
\end{equation} 
By lemma $A.1$ \cite{bri1}, the groupoid of $T$-valued points of $Z$
can be described as follows. The objects are short exact sequences of
$T$-flat sheaves on $T\times Y$ of the form
$$
0\to A \to G \to B \to 0
$$
such that  $A$ and $B$ define families of sheaves on $Y$ lying in the
subcategories $\Pp$ and $\Qp$ respectively. The morphisms are
isomorphisms of short exact sequences. The composition,
$$
g = b\circ f: Z \to \MM
$$
sends a short exact sequence to the object $G$. Since $\Pp$ and $\Qp$ are subcategories of $\CCC$, it follows that the composition factors through $\CCC$. This morphism induces an equivalence on $\CC$-valued points because of the torsion pair property: every object $G$ of $\CCC$ fits into a unique short exact sequence of the form (\ref{ugh}). Thus, the identity follows from the relations in the infinite Hall algebra.
\end{pf}

We will need a framed version of the previous lemma.

\begin{lemma} \label{frame}
$$
\oneo{\CCC} = \oneo{\Pp} \ast \oneo{\Qp} .
$$
\end{lemma}

\begin{pf}
Form the following Cartesian diagram:
\[
\xymatrix{ 
U \ar[r]^p \ar[d] & V \ar[r]^j \ar[d] & \MM^{(2)} \ar[r]^b \ar[d]^{(a_1,a_2)} & \MM \\
\Pp (\OO) \times\Qp(\OO) \ar[r]^{q\times \text{id}}  &\Pp \times\Qp(\OO) \ar[r] &\MM\times\MM .
}
\]

Then $\onePpO \ast \oneQpO$ is represented by the composite map
$b\circ j \circ p: U \to \MM$. Since $R^{1}\pi _{*}P=0$ for all $P\in
\Pp $, the argument of lemma 2.5 in \cite{bri2} implies that the map
$\Pp (\OO) \to \Pp$ is a Zariski fibration with fibre $H^0(P)$ over a
point $P \in \Pp$. By pullback, the same is true of the morphism $p$.

The groupoid of $T$-valued points of $V$ can be described as follows. The objects are short exact sequences of $T$-flat sheaves on $T\times Y$ of the form
$$
0 \to P \to G \to B \to 0
$$
such that $P$ and $B$ define flat families of objects in $\Pp$ and $\Qp$ respectively, together with a map $\OO_{{T\times Y}} \to B$.  We represent the objects of $V$ as diagrams of the form:
\[
\begin{CD}
  & @.   &   & @.   &   & @.   &\OO_{T\times Y} & @.   & \\
  & @.   &   & @.   &   & @.   & @VVV           & @.   & \\
0 & @>>> & P & @>>> & G & @>>> & B              & @>>> & 0.
\end{CD}
\]

Consider the stack $Z$ from lemma~\ref{basiclemma} with its map $Z \to \MM$. Form the diagram
\[
\begin{CD}
W & @>h>> & \MM  (\OO) \\
@VVV & & @ VqVV \\
Z & @>g>> & \MM .
\end{CD}
\]
Since $g$ induces an equivalence on $\CC$-valued points, so does $h$, so that the element $\oneo{\CCC}$ can be represented by the map $q\circ h$.

We represent the objects of $W$ as diagrams of the form:
\[
\begin{CD}
  & @.     &   & @.     &  \OO_{T\times Y}& @.         &   & @.   & \\
  & @.     &   & @.     & @VV\delta V     & @.         &   & @.   & \\
0 & @>>>   & P & @>>>   & G               & @>\beta >> & B & @>>> & 0.
\end{CD}
\]
Setting $\gamma = \beta \circ \delta$ defines a map of stacks $W \to V$, which is a Zariski fibration with fibre an affine model of the vector space $H^0(P)$ over a sheaf $P$. 

Now, since $H^1(P)=0$, $U \to V$ is a Zariski fibration with fibre
$H^0(P)$, hence they represent the same element of the Hall algebra,
namely $\one{\CCC}^{\OO} $.
\end{pf}

\begin{lemma} \label{hilb}
$$
\one{\CCC}^{\OO} = \HHH \ast \one{\CCC}
$$
\end{lemma}
This is lemma 4.3 of \cite{bri2}. Intuitively, this amounts to the fact that every map $\OO \stackrel{\gamma}{\to} G$ factors uniquely into a surjection $\OO \to \im (\gamma)$ and an inclusion $\im (\gamma) \to G$. The following lemma is a restriction of the previous to the substack $\Pp$.

\begin{lemma} \label{stepthree}
$$
\onePpO = \HHe \ast \onePp
$$
\end{lemma}

\begin{pf}
Form the Cartesian diagram:
\[
\begin{CD}
Z & @>>> &\MM^{(2)}&@>b>>&\MM \\
 @VVV & &@VV(a_1,a_2)V \\
\HHe \times \Pp &@>i>> &\MM\times\MM.
\end{CD}
\] 
The groupoid of $T$-valued points of $Z$ may be described as
follows. The objects are short exact sequences of $T$-flat sheaves of
$T\times Y$
$$
0\to A \to G \to B \to 0
$$
such that for all geometric points $t\in T$, $B_{t} \in \Pp$, and
$A_{t}$ is supported on exceptional fibres, together with an
epimorphism $\OO_Y \to A_{t}$. We can represent these objects as
diagrams of the form
\begin{equation*} 
 \xymatrix@C=1.5em{  &\OO_{T \times Y} \ar[d]^{\gamma}
\\
0 \ar[r] &A \ar[r]^{\alpha}& G  \ar[r]^\beta &B \ar[r] &0. }
\end{equation*}

Let $t\in T$ be an arbitrary geometric point.  Since $\OO_Y \to A_{t}$
is an epimorphism, we know that $A_{t}$ is of the form $\OO_{C_{t}}$
for some one-dimensional subscheme $C_{t}$ of $Y$. By the proof of
proposition~\ref{prop1}, we know that $R\pi_*\OO_{C_{t}}=0$. Since
$A_{t}=\OO_{C_{t}}$ has exceptional support, it follows that
$\pi_*A_{t}$ is a zero-dimensional sheaf, hence $A_{t}\in \Pp$. Since
$B_{t}\in \Pp$ by design, and since $\Pp$ is closed under taking
extensions, we conclude that in any such short exact sequence, $G_{t}\in
\Pp$. There is a map $h:Z \to \Pp (\OO) $ sending the above diagram to
the composite map
$$
\delta = \alpha \circ \gamma : \OO_{T\times Y} \to G.
$$
This morphism $h$ fits into a commuting diagram of stacks
$$ 
\xymatrix@C=1em{
Z\ar[rr]^{h} \ar[dr]_{b\circ f} &  & \Pp (\OO) \ar[dl]^{q} \\  &\MM }
$$
We argue that the map $h$ then induces an equivalence on $\CC$-valued
points. Suppose $\OO_{T\times Y} \stackrel{\delta}{\to} G$ is an
arbitrary map of sheaves, with $G$ defining a family of sheaves in
$\Pp$. Then we get the following diagram:

\begin{equation*} 
 \xymatrix@C=1.5em{  &\OO_{T \times Y} \ar[d]^{\gamma}
\\
0 \ar[r] &\im (\delta) \ar[r]& G  \ar[r] &\cok (\delta) \ar[r] &0. }
\end{equation*}

Since $G_{t}\in \Pp$ we know that the one-dimensional component of its
support is exceptional, hence $\im (\delta)|_{t}$ is also exceptional,
so that $\OO_{T\times Y} \to \im (\delta)$ defines a family of objects
in $\Hilbe$. As well, we know that $\Pp$ is closed under taking
quotients, so $\cok (\delta)$ is in $\Pp$. This completes the proof.
\end{pf}

Morally, the next lemma is similar to lemma~\ref{hilb} since $\HHp$ may be thought of as the surjections $\OO_Y \to G$ in a tilt of the abelian category generated by $\OO$ and $\CCC$. We provide a direct proof since we have not constructed this tilt.
\begin{lemma} \label{oneqhponeq} \label{stepone}
$$
\oneQpO = \HHp \ast \oneQp.
$$
\end{lemma}

\begin{pf}
Form the following Cartesian diagram:
\[
\begin{CD}
Z & @>f>> &\MM^{(2)}&@>b>>&\MM \\
 @VVV & &@VV(a_1,a_2)V \\
\HHp\times\Qp &@>i>> &\MM\times\MM
\end{CD}
\] 
The groupoid of $T$-valued points of $Z$ is described as follows. The objects are short exact sequences of $T$-flat sheaves on $T\times Y$
$$
0 \to A \to G \to B \to 0
$$
with the property that $B \in \Qp$, together with a map $\OO_{T\times
Y} \to A$ that pulls back to a \ps\hbox{ }pair $\OO_Y \to A_t$ for
every $t\in T$. We can represent these objects as diagrams of the
form:

\begin{equation} \label{obob}
 \xymatrix@C=1.5em{  &\OO_{T \times Y} \ar[d]^{\gamma}
\\
0 \ar[r] &A \ar[r]^{\alpha}& G  \ar[r]^\beta &B \ar[r] &0. }
\end{equation}
Since $A$ and $B$ are objects of $\Qp$, and $\Qp$ is closed under extensions, we conclude $G\in \Qp$.
Thus, there is a map $h:Z \to \Qp (\OO)$ sending the above diagram to the composite map 
$$
\sigma  = \alpha \circ \gamma : \OO_{T\times Y} \to G.
$$ 
 This map $h$ fits into a commuting diagram of stacks

\begin{equation} \xymatrix@C=1em{
Z\ar[rr]^{h} \ar[dr]_{b\circ f} &  & \Qp (\OO) \ar[dl]^{q} \\  &\MM }
\end{equation}
 The map $h$ then induces an equivalence on $\CC$-valued points because
of the following argument. Let $\OO_Y \stackrel{\sigma}{\to} G$ be an
arbitrary map, with $G\in \Qp$. We need to produce a diagram,
$$
\xymatrix@C=1.5em{  &\OO_{Y} \ar[d]^{\gamma}
\\
0 \ar[r] &A \ar[r]^{\alpha}& G  \ar[r]^\beta &B \ar[r] &0, }
$$
with $\OO_Y \to A$ a $\pi$-stable pair, $B\in\Qp$, and $\alpha \circ \gamma = \sigma$. 

Consider the cokernel $K$ of $\sigma$. Since $(\Pp, \Qp)$ is a torsion pair, we know that $K$ fits into a short exact sequence
$$
0 \to P \to K \to Q \to 0
$$
where $P\in\Pp$ and $Q\in\Qp$. Let $G \stackrel{c}{\to} K$ be the canonical map from $G$ to the cokernel $K$ of $\sigma$, and let $G\stackrel{d}{\to} Q$ be the composition of $G\to K$ and $K \to Q$. Define $A$ to be the kernel of $d$. Consider the following diagram.
$$
\xymatrix@C=1.5em{
0\ar[r] & \OO_Y \ar @{.>}[d] \ar[r]^{=} & \OO_Y \ar[r] \ar[d]^{\sigma} & 0\ar[d]\\
0\ar[r] & A \ar @{.>}[d] \ar[r]  & G \ar[r]^{d} \ar[d]^{c} & Q \ar[d]^{=} \ar[r] & 0 \\
0\ar[r] & P \ar @{.>}[d] \ar[r] & K\ar[r]\ar[d] & Q\ar[d] \ar[r] & 0\\
 &0 & 0& 0. &
}
$$
We know that $Q \in \Qp$, and a diagram chase proves that the dotted
vertical morphisms exist and that $P$ is the cokernel of $\OO_Y \to
A$. The sheaf $A$ is a subsheaf of $G\in \Qp$, and $\Qp$ is closed
under taking subsheaves, so $A \in \Qp$. This proves that $\OO_Y \to
A$ is a $\pi$-stable pair, and thus $h$ is a surjection on $\CC
$-valued points. Moreover, the above diagram is uniquely determined up
to isomorphism (since the exact sequence $P\to K\to Q$ is unique up to
isomorphism) and consequently the preimage of $\OO _{Y}\to G$ under
$h$ is unique up to isomorphism. Thus $h$ is a geometric bijection and
thus a (constructible) equivalence of stacks \cite[Lemma~3.2]{bri1}. 
\end{pf}


We end this section by giving a fake proof of
theorem~\ref{therecanbeonlyone} that depends on a fake integration
map. In fact, no such integration map is known to exist, but if there
was one, the proof of our theorem would be simpler. As it stands, we have
a chapter dedicated to convergence issues to get around the fact that
no such integration map exists on the infinite type Hall algebra. It
is our hope that this fake proof will make the true one easier to
follow.

\begin{fpf}  
Lemma 4.3 of \cite{bri2} proves
$$
\Hleq \ast \one{\CCC} =\oneo{\CCC}.
$$
Using lemma~\ref{basiclemma}, we may rewrite $\Hleq \ast \one{\CCC}$:
$$
\Hleq \ast \one{\CCC} = \Hleq \ast \one{\Pp} \ast \one{\Qp}.
$$
Lemma~\ref{frame} allows us to write
$$
\oneo{\CCC} = \oneo{\Pp} \ast \oneo{\Qp}.
$$
Putting these together yields
$$
\Hleq \ast \one{\Pp} \ast \one{\Qp} = \oneo{\Pp} \ast \oneo{\Qp}.
$$
Applying lemma~\ref{stepthree}
$$
\Hleq \ast \one{\Pp} \ast \one{\Qp} = \HHe \ast \one{\Pp} \ast \oneo{\Qp}.
$$

From lemma~\ref{stepone}, we have
\begin{equation} 
 \oneQpO=\HHp \ast \oneQp .
\end{equation}
so we get
$$
\Hleq \ast \one{\Pp} \ast \one{\Qp} = \HHe \ast \one{\Pp} \ast \HHp \ast \oneQp.
$$
Now, for reasons we will explain in the next section, $\one{\Pp}$ and $\one{\Qp}$ are invertible in the Hall algebra. We may therefore cancel the copies of $\one{\Qp}$ and isolate $\Hleq$.
$$
\Hleq  = \HHe \ast \one{\Pp} \ast \HHp \ast \one{\Pp}^{-1}.
$$

The elements $\Hleq, \HHe, \HHp$ all lie in the subalgebra $\HRA$
since they are represented by (constructible) schemes. As we will see
in the next section, conjugation by $\one{\Pp}$ induces a Poisson
homomorphism of $\HRA$ of the form: identity + terms expressed in the
Poisson bracket. Since the Poisson bracket of the polynomial ring is
trivial, these terms vanish when we apply the fake integration map,
and we are left with
$$
I (\HH) = I (\HHe) \cdot I (\HHp)
$$
or equivalently,
$$
\frac{I(\HH)}{I(\HHe)} = I(\HHp).
$$

Up to signs arising from lemma~\ref{bflemma}, the ``polynomials"
$I(\HH)$, $I(\HHe)$, and $I(\HHp)$ are the generating series of
$\DT(Y)$, $\DT_{exc}(Y)$, and $\TP (Y)$, respectively, and we see that
the above equation is the formula claimed in
theorem~\ref{therecanbeonlyone}. \end{fpf}

The true proof will follow precisely these steps, fully justified, and with the appropriate convergence arguments. The next chapter describes the Laurent Hall algebra, which does have an integration map.

\section{Equations in the Laurent Hall algebra and the true proof}

\subsection*{Laurent subsets}
In this section, we will formally modify the algebra $\HA$ and its integration map so that the modified integration map on the modified algebra takes valued in power series. This section is a summary of sections 5.2 and 5.3 of \cite{bri2}.

\begin{definition} A subset $S \sub \Delta$ is \emph{Laurent} if for all $\beta \in N_1(Y)$, the collection of elements of the form $(\beta, n) \in S$ is such that $n$ is bounded below. 
\end{definition}
 Let $\Phi$ denote the set of all Laurent subsets. It has the following properties:
\begin{enumerate}
\item if $S, T \in \Phi$ then so it $S + T = \{ \alpha +\beta : \alpha \in S, \beta \in T\}$
\item if $S, T \in \Phi$, and $\alpha \in \Delta$ then there are only finitely many ways to write $\alpha = \beta + \gamma$ such that $\beta \in S$ and $\gamma \in T$.
\end{enumerate}

Given a ring $A$ graded by $\Delta$, $A = \bigoplus_{\gamma \in \Delta}A_\gamma$, we can use the Laurent subset to define a new algebra, which we will denote $A_\Phi$. Elements of $A_\Phi$ are of the form
$$
a=\sum_{\gamma \in S}a_\gamma
$$
where $S \in \Phi$, and $a_\gamma \in A_\gamma \sub A$. Given an element $a\in A_\Phi$ as above, we define $\rho_\gamma (a) = a_\gamma \in A$. (Here, our notation differs from \cite{bri2}, since we are using the symbol $\pi$ for the map $\pi: Y \to X$.) The projection operator $\rho$ allows us to define a product $\ast$ on $A_\Phi$ by:
$$
\rho_\gamma(a\ast b) = \sum_{\gamma = \alpha + \beta} \rho_\alpha (a)\ast \rho_\beta(b).
$$
$A_\Phi$ admits a natural topology that may be identified by declaring a sequence $(a_j)_{j\in \NNN} \sub A_\Phi$ to be convergent if for any $(\beta, n) \in \Delta$, there exists an integer $K$ such that for all $m<n$ 
$$
i, j > K \Rightarrow \rho_{(\beta, m)}(a_i) = \rho_{(\beta, m)}(a_j).
$$

\begin{lemma} If A is a $\CC$-algebra and $a\in A_\Phi$ satisfies $\rho_0(a) =0$ then any series
$$
\sum_{j\geq 1}c_ja^j
$$
with coefficients $c_j \in \CC$ is convergent in the topological ring $A_\Phi$.
\end{lemma}
See Lemma 5.3 of \cite{bri2} for a proof.

Given two $\Delta$-graded algebras $A$ and $B$, and a morphism $f: A \to B$ that preserves the $\Delta$-grading, we get an induced continuous map 
$$
f_\Phi : A_\Phi \to B_\Phi
$$
by defining
$$
\rho_\gamma (f_\Phi (a)) = f(\rho_\gamma (a)).
$$

Applying this process to the map of $\Delta$-graded algebras $I: \HSC \to \CC [\Delta]$ yields a continuous map 
$I_\Phi: \HSC_\Phi \to \CC [\Delta]_\Phi$. We call $\HSC_\Phi$ the \emph{Laurent} Hall algebra; it, too, is equipped with a Poisson bracket. 

\begin{definition} A morphism of stacks $f: W \to \CCC$ is $\Phi$-\emph{finite} if
\begin{enumerate}
\item[a)] $W_\alpha = f^{-1}( \CCC_\alpha)$ is of finite type for all $\alpha \in \Delta$, and
\item[b)]there is a Laurent subset $S \sub \Delta$ such that $W_\alpha$ is empty unless $\alpha \in S$.
\end{enumerate}
\end{definition}

A $\Phi$-finite morphism of stacks $f: W \to \CCC$ defines an element of $\HSC_\Phi$ by the formal sum
$$
\sum_{\alpha \in S}[W_\alpha \stackrel{f}{\to} \CCC].
$$

In \cite{bri2}, it is shown that $[\Hilb \to \MM]$ is $\Phi$-finite, and hence defines an element in $\HSC_\Phi$.

\begin{lemma} \label{integrate}
The maps 
$$
\Hilbp \to \CCC, \hspace{0.6cm}  \Hilb \to \CCC, \hspace{0.6cm}  \Hilb_{exc} \to \CCC,
$$
are $\Phi$-finite. The corresponding elements $\HHp, \HHH$ and $\HHe$ of $\HSC_\Phi$ satisfy
$$
I_\Phi (\HHp) = \sum_{(\beta,n)\in \Delta} (-1)^n \TP(\beta, n)x^\beta q^n = \TP (Y)(x, -q),
$$
where we have written $x^\beta = q^{(\beta, 0)}$ and $q=q^{(0,1)}.$ Similarly, 
$$
I_\Phi (\HHH) = \sum_{(\beta,n)\in \Delta} (-1)^n \DT(\beta, n)x^\beta q^n = \DT (Y) (x, -q)
$$
$$
I_\Phi (\HHe) = \sum_{(\beta,n)\in \Delta, \pi_* \beta=0}
(-1)^n DT(\beta, n)x^\beta q^n = \DT_{exc}(Y) (x, -q).
$$
\end{lemma} 

\begin{pf} Lemma~\ref{ftalgspace} proves that $\Hilbp$ is
constructible and locally of finite type, and is of finite type once
the Chern character is fixed.  As well, the set of elements $\alpha
\in \Delta$ for which $\Hilbpa$ is non-empty is Laurent. To prove
this, it suffices to show that for any curve class $\beta$, there
exists an integer $N$ such that for any $n<N$, the moduli space
$\Hilbpbn$ is empty. Fix a curve class $\beta$, and consider all \ps
pairs $\OO_Y \to G$ in that class. There is the associated short exact
sequence,
$$
0 \to \OO_C \to G \to P \to 0
$$
and since there are only finitely many decompositions $\beta =\beta_1 + \beta_2$ of $\beta$ into a sum of effective curve classes, we lose no generality in fixing the curve class of $\OO_C$ and $P$. Now the structure sheaf $\OO_C$ lives in a Hilbert scheme, and the set of elements $(\beta_1, n_1) \in \Delta$ for which $\Hilb^{(\beta_1, n_1)}$ is non-empty is Laurent, so there is a ``minimal" Euler characteristic of $\OO_C$, which we denote by $N_1$. As for $P$, the Leray spectral sequence shows that $\chi (P)\geq 0$ (see lemma~\ref{bound}). This proves that we may take $N=N_1$, and for any $n<N$, $\Hilbpbn$ is empty.

The formulae then follow from lemma~\ref{bflemma} and Behrend's description of DT invariants as a weighted Euler characteristic. 

$\Hilb_{exc}$ is a subscheme of $\Hilb$, so the desired properties follow from \cite[lemma 5.5]{bri2}.
\end{pf}

\begin{lemma} Let $I \subset (-\oo, +\oo] \times (-\oo, +\oo]$ be an interval bounded from below. Then
$$
\one{\SS (I)} \to \CCC
$$
is $\Phi$-finite.
\end{lemma}
\begin{pf} 
Since this holds for Gieseker stability (\cite[theorem 3.3.7]{hule}), it suffices to prove that for any $b=(b_1, b_2) \in (-\oo, +\oo] \times (-\oo, +\oo]$ there exists a number $M$ such that the family of all $G$ with $\mup (G) \geq b$ satisfies $\mu (G) \geq M$. Here, $\mu$ stands for Gieseker slope stability, namely
$$
\mu (G) = \frac{\chi (G)}{\beta \cdot L}.
$$
\begin{itemize}
\item[Case 1: ] $\chi(G)>0$

Here, we have
$$
0< \frac{\chi(G)}{\beta \cdot L} = \mu (G),
$$
so we may take $M=0$ in this case.
\item[Case 2: ] $\chi(G)<0$

Now, since $\beta \cdot \wt{H} \leq \beta \cdot L$
(c.f. definition~\ref{defn: mupi stability}) and $\chi(G)<0$, we have
$$
b_1 \leq \frac{\chi(G)}{\beta \cdot \wt{H} } \leq \frac{\chi(G)}{\beta \cdot L}.
$$
In this case, we may take $M=b_1$.
\item[Case 3: ] $\chi(G) =0$

In this case, $\mu(G)$ is either $0$ or $+\oo$, so we may take $M=0$ in this case. 
\end{itemize}
\noindent Case-by-case analysis reveals that we may use $M = \min \{ 0, b_1\}$.
\end{pf}

\subsection*{Equations in the Laurent Hall algebra} In this section, following \cite{bri2} we establish equations in $\HSC_\Phi$, and ultimately prove theorem~\ref{therecanbeonlyone}. 

\begin{lemma} Let $\mu \in (-\oo, +\oo] \times (-\oo, +\oo]$ such that $\mu < \twoo$. Then the following equality holds in $\HSC_\Phi$:
$$
\one{ \smloo} = \one{\Pp} \ast \one{\smltwoo}.
$$
\end{lemma}

\begin{pf}
Form the following Cartesian diagram:
$$
\xymatrix{
Z \ar[r]^f \ar[d] & \MM^{(2)} \ar[r]^c \ar[d]^{(l,r)} & \MM \\
\Pp\times \smltwoo \ar[r] & \MM\times\MM &
}
$$ 
$T$-valued points of $Z$ are short exact sequences $0\to A \to G \to B
\to 0$ of $T$-flat sheaves on $T\times Y$ such that $A$ defines a
family of objects in $ \Pp$ and $B$ a family in $\smltwoo$. By
lemma~\ref{pqss}, we know $\twoo \leq \mup (A) \leq \oo$. Now by
\cite[lemma~6.2]{bri2}, we know that $G$ defines a family of objects
in $\smloo.$

Now let $G \in \smloo$. If $G\in \Pp$ or $\smltwoo$ then we are done, since then $G$ will be an extension where one term is zero (recall that all $\SS (a<\square <b)$ include the zero objects). Otherwise, let
$$
0=G_0 \sub G_1 \sub \ldots \sub G_{N-1} \sub G_N=G
$$
be the Harder-Narasimhan filtration of $G$, let $Q_i$ be the associated quotients. Then there exists an index $j\in \NNN, 1<j<N$ such that for all $i<j, \mup(Q_i)\geq \twoo$ and $\mup(Q_j)<\twoo$. Finally, by the uniqueness of the Harder-Narasimhan filtration, we have that $G_j\in \Pp$ and $G/G_j\in \smltwoo$. 

\end{pf}

\begin{remark} The proof of this lemma is strikingly similar to the proof of Lemma~\ref{basiclemma}. This is no coincidence. The above is actually just a minor refinement of Lemma~\ref{basiclemma} which says that we may cut off the tail end of $\Qp$ and have the corresponding result still hold. As we go on, we will be less explicit about the proofs of lemmas when the argument has been already made in the infinite-type case. 
\end{remark}

\begin{lemma} Let $\mu \in (-\oo, +\oo]\times (-\oo, +\oo].$ Then, as $\mu \to -\oo$, we have
$$
\HHH \ast \one{\smloo} - \oneo{\smloo} \to 0.
$$
\end{lemma}
\begin{pf}
Fix $(\beta, n) \in \Delta$. Then there are only finitely many
decompositions $(\beta, n) = (\beta_1, n_1) + (\beta_2, n_2)$ such
that both $\rho_{(\beta_1, n_1)}(\HHH)$ and $\rho_{(\beta_2, n_2)}(
\one{\smloo} )$ are non-zero. This follows from the fact that there
are only finitely many decompositions $\beta = \beta_1 +\beta_2$ with
both $\beta_i$ effective. Now for each fixed $\beta$, there exist
finitely many $n$ such that both $\rho_{(\beta_1, n_1)}(\HHH)$ and
$\rho_{(\beta_2, n_2)}( \one{\smloo} )$ are non-zero).

By the boundedness of the Hilbert scheme, we may assume that $\mu$ is small enough so that for any of the decompositions, $\beta = \beta_1 + \beta_2$, all points $\OO_Y \to A$ of $\Hilb^{(\beta_1, n_1)}$ satisfy $A\in \SSS ( \mu \leq \square\leq \oo)$. Consider a diagram of sheaves, 
 \begin{equation}
 \xymatrix@C=1.5em{  &\OO_{Y} \ar[d]^{\gamma}
\\
0 \ar[r] &A \ar[r]^{\alpha}& G  \ar[r]^\beta &B \ar[r] &0 }
\end{equation}
with $\OO_Y \to A$ in $\Hilb^{(\beta_1, n_1)}$ and $\text{ch}([G]) = (\beta, n)$.

Now $G \in \SSS ( \mu \leq \square \leq \oo)$ if and only if $B \in \SSS ( \mu \leq \square\leq \oo)$. Since Bridgeland proves \cite[prop 6.5]{bri2}
$$
\rho_{(\beta, n)} (\HHH \ast \one{\smloo} ) = \rho_{(\beta, n)} (\oneo{\smloo} ),
$$
the claim is proven. 
\end{pf}

\begin{lemma}  Let $\mu \in (-\oo, +\oo]\times (-\oo, +\oo).$ Then, as $\mu \to -\oo$, we have
$$
\HHp \ast \one{\smltwoo} - \oneo{\smltwoo} \to 0.
$$
\end{lemma}
\begin{pf}
 Fix $(\beta, n) \in \Delta$. Then there are only finitely many
decompositions $(\beta, n) = (\beta_1, n_1) + (\beta_2, n_2)$ such
that both $\rho_{(\beta_1, n_1)}(\HHp)$ and $\rho_{(\beta_2, n_2)}(
\one{\smltwoo} )$ are non-zero.

By the boundedness of the moduli space of $\pi$-stable pairs, we may assume that $\mu$ is small enough so that for any decompositions, $\beta = \beta_1 + \beta_2$, all points $\OO_Y \to A$ of $\Hilbp^{(\beta_1, n_1)}$ satisfy $A\in \SSS ( \mu \leq \square < \twoo)$.
Consider a diagram of sheaves, 
 \begin{equation}
 \xymatrix@C=1.5em{  &\OO_{Y} \ar[d]^{\gamma}
\\
0 \ar[r] &A \ar[r]^{\alpha}& G  \ar[r]^\beta &B \ar[r] &0 }
\end{equation}
with $\OO_Y \to A$ in $\Hilbp (\beta_1, n_1)$ and $[G] = (\beta, n)$. Using that $\SSS (I) \ast \SSS (I) \subset \SSS (I)$, we see that $B \in \SSS ( \mu \leq \square< \twoo)$ if and only if $G \in \SSS ( \mu \leq \square< \twoo)$.  Now since $A \in \Qp$ and $B \in \Qp$, we have that $G \in \Qp$. Composing the map $\OO_Y \to A$ with the map $A \to G$, yields a map $\OO_Y \to G$; this represents an object of $\oneo{\Qp}$. The proof of lemma~\ref{oneqhponeq}, that
$$
\oneo{\Qp} = \HHp \ast \one{\Qp},
$$
can be easily adapted to now prove that 
$$
\rho_{(\beta, n)}(\HHp \ast \one{\smltwoo} ) = \rho_{(\beta, n)} ( \oneo{\smltwoo} )
$$
This completes the proof that 
$$
\HHp \ast \one{\smltwoo} - \oneo{\smltwoo} \to 0
$$
as $\mu \to -\infty$.
\end{pf}

\begin{prop} We have the following equality in the Laurent Hall algebra, $\HSC_\Phi$:
$$
\HHH \ast \one{\Pp} = \HHe \ast \one{\Pp} \ast \HHp.
$$
\end{prop}
\begin{pf}
Using  $\one{\smloo} = \one{\Pp} \ast \one{\smltwoo}$ and $\one{\smloo}^{\OO} = \onePponeO \ast \one{\smltwoo}^{\OO},$ we can rewrite 
$$
\HHH \ast \one{\smloo} - \one{\smloo}^{\OO} \to 0
$$
as
$$
\HHH \ast \one{\Pp} \ast \one{\smltwoo} - \oneo{\Pp} \ast \oneo{\smltwoo} \to 0,
$$
as $\mu \to -\infty$. 

Multiplying $\HHp \ast \one{\smltwoo} - \oneo{\smltwoo} \to 0$ on the left by $\oneo{\Pp}$, and rewriting using $\oneo{\Pp} = \HHe \ast \one{\Pp}$ yields
$$
\HHe \ast \one{\Pp} \ast \HHp \ast \one{\smltwoo} - \oneo{\Pp} \ast \oneo{\smltwoo} \to 0
$$
as $\mu \to -\infty$. Hence
$$
\HHe \ast \one{\Pp} \ast \HHp \ast \one{\smltwoo} - \HH \ast \one{\Pp} \ast \one{\smltwoo} \to 0
$$
as $\mu \to -\infty$. Since $\one{\smltwoo}$ is invertible, we can cancel it from both sides:
$$
\HHH \ast \one{\Pp} = \HHe \ast \one{\Pp} \ast \HHp.
$$
\end{pf}

\subsection*{The proof of theorem~\ref{therecanbeonlyone}}
We first collect results. The next proposition is theorem 3.11 of \cite{joso}, and is a very deep result whose proof depends on all the full power of the formalism of \cite{joy1, joy2, joy3, joy4, joy5}.

\begin{prop} \label{epsilon}
For each slope $\mu \in ( (-\oo, -\oo), (+\oo, +\oo)]$, we can write 
$$
\one{\SSS (\mu)} = \exp (\epsilon_\mu) \in \HC_\Phi
$$
with $\nu_\mu = [\CC^*]\cdot \epsilon_\mu \in \HRC_\Phi$ a regular element.
\end{prop}
\begin{pf} The proof is identical to that of \cite[theorem 6.3]{bri2}. Bridgeland uses Joyce's machinery, which applies in our case just as it does in his.
\end{pf}
The following corollary corresponds to Bridgeland's 6.4 \cite{bri2}.
\begin{corollary} \label{poissauto}
For any $\mu \in ( (-\oo, -\oo), (+\oo, +\oo)]$, the element $\one{\SSS (\mu)} \in \HC_\Phi$ is invertible, and the automorphism
$$
\Ad_{\one{\SSS (\mu)}} : \HC_\Phi \to \HC_\Phi
$$
preserves the subring of regular elements. The induced Poisson automorphism of $\HSC$ is given by
$$
\Ad_{\one{\SSS (\mu)}} = \exp\{\eta_\mu, -\}.
$$
\end{corollary}
\begin{pf} The proof of this is identical to that of corollary 6.4 of \cite{bri2}.
\end{pf}

Now we can prove theorem~\ref{therecanbeonlyone}. We have
$$
\HHH \ast \one{\Pp} = \HHe \ast \one{\Pp} \ast \HHp.
$$
Rearranging yields
$$
\HHH = \HHe \ast \one{\Pp} \ast \HHp \ast (\one{\Pp})^{-1}.
$$
By lemma~\ref{pqss}, we can write $\one{\Pp} = \SSS ( \twoo \leq \square \leq \oo)$, and by lemma~6.2 of \cite{bri2}, we can write 
$$
\SSS ( \twoo \leq \square \leq \oo) = \prod_{\twoo \leq \mu \leq \oo} \one{ \SSS ( \mu )}.
$$
In \cite[lemma 6.2]{bri2}, it is explained that given an interval $J \sub (-\infty, +\infty] \times (-\infty, +\infty]$ that is bounded below, and an increasing sequence of finite subsets 
$$
V_1 \sub V_2 \sub \ldots \sub J
$$
the sequence $\one{\SSS (V_j)}$ converges to $\one{\SSS (J)}$, where $\one{\SSS (V_j)}$ is defined to be 
$$
\prod_{v\in V_j} \one{\SSS (v)},
$$
where the product is taken in descending order of slope. So, letting $J$ denote the interval of slopes between $\twoo$ and $\oo$, including $\twoo$ and excluding $\oo$, we can write
$$
\HHH = 
$$
$$\HHe \ast \lim_{\text{finite } V \sub J} \one{\SSS (\mu_N)} \ast \ldots \ast \one{\SSS (\mu_1)} \ast \HHp \ast (\one{\SSS (\mu_1)})^{-1} \ast \cdots \ast (\one{\SSS (\mu_N)})^{-1},
$$
where $\mu_i$ enumerate all the elements of $V$. Using proposition~\label{epsilon} and corollary~\label{poissauto}, we can rewrite
$$
\HHH = 
$$
$$\HHe \ast \lim_{\text{finite } V \sub J} \exp(\{\eta_{ \mu_N }, \exp\{\eta_{\mu_{N-1}}, \ldots \exp \{\eta_{\mu_1}, -\} \ldots \} (\HHp).
$$
Now hitting this equation with the integration map yields
$$
I_\Phi (\HHH) = 
$$
$$
I_\Phi (\HHe) \cdot I_\Phi \left( \lim_{\text{finite } V \sub J} \exp(\{\eta_{\mu_N}, \exp\{\eta_{\mu_{N-1}}, \ldots \exp \{\eta_{\mu_1}, -\} \ldots \} (\HHp) \right).
$$
The integration map commutes with limits since it is continuous, thus
$$
I_\Phi (\HHH) =
$$
$$
I_\Phi (\HHe) \cdot \lim_{\text{finite } V \sub J} I_\Phi \left( \exp (\{ \eta_{\mu_N}, \exp( \{ \eta_{\mu_{N-1}}, \ldots \exp( \{ \eta_{\mu_1}, -\}) \ldots \}) (\HHp) \right).
$$
Now, the Poisson bracket is a commutator which is trivial in the ring of Laurent series, so it vanishes after applying the integration map, and we are left with
$$
I_\Phi (\HHH) = I_\Phi (\HHe) \cdot I_\Phi (\HHp).
$$
Applying lemma~\ref{integrate}, we get
$$
\DT(Y)(x, -q) = \DT_{exc}(Y)(x, -q) \cdot \TP (Y)(x, -q)
$$
and substituting $q$ for $-q$ yields
$$
\DT(Y) = \DT_{exc}(Y) \cdot \TP (Y),
$$
which is what we set out to prove. \hfill $\blacksquare$ 



\begin{thebibliography}{00}

\bibitem{ba} A. Bayer: \emph{ Polynomial Bridgeland stability conditions and the large volume limit}, arXiv: math.AG/0712.1083.

\bibitem{ba2} A. Bayer: (\emph{In preparation})

\bibitem{beh} K. Behrend: \emph{Donaldson-Thomas type invariants via microlocal geometry}, Ann. of Math. (2) 170 (2009), no. 3, 1307–1338, arXiv: math.AG/0507523.

\bibitem{bebr} K. Behrend, and J. Bryan: \emph{Super-rigid Donaldson-Thomas invariants}, Mathematical Research Letters, 14(4):559--571, 2007, arXiv: math.AG/0601203.

\bibitem{bf1} K. Behrend, and B. Fantechi: \emph{The intrinsic normal cone}, Invent. Math. 128 (1997), 45-88.

\bibitem{bf2} K. Behrend, and B. Fantechi: \emph{Symmetric obstruction theories and Hilbert schemes of points on threefolds,} Algebra \& Number Theory 2 (2008), no. 3, 313--345, arXiv: math.AG/0512556. 

\bibitem{bbd} A. Be\u{\i}linson, J. Bernstein, P. Deligne: \emph{Faisceaux perverse}, Analysis and topology on singular spaces, I (Luminy, 1981), 5–171, Ast\'{e}isque, 100, Soc. Math. France, Paris, 1982.

\bibitem{bosa} S. Boissiere, and A. Sarti; \emph{Contraction of excess fibres between the McKay correspondences in dimensions two and three,} Ann. Inst. Fourier, 57(6):1839--1861,2007 arXiv: math.AG/0504360.

\bibitem{briflop} T. Bridgeland: \emph{Flops and derived categories}, Invent. Math. 147, 613–632 (2002). math.AG/0009053.

\bibitem{brs} T. Bridgeland: \emph{Stability conditions on triangulated categories,} Ann. of Math, Vol 166, 317-345, 2007.

\bibitem{bri1} T. Bridgeland: \emph{An introduction to motivic Hall algebras}, preprint arXiv: math.AG/1002.4372v1. 

\bibitem{bri2} T. Bridgeland: \emph{Hall algebras and curve-counting invariants,} J. Amer. Math. Soc. 24 (2011), no. 4, 969–998.

\bibitem{bkr} T. Bridgeland, A. King, and M. Reid: \emph{The McKay correspondence as an equivalence of derived categories,} J. Amer. Math. Soc. 14 (2001), no. 3, 535–554 (electronic). 

\bibitem{bcy} J. Bryan, C. Cadman, and B. Young: \emph{The orbifold topological vertex}, Advances in Mathematics, 229 (1), pg. 531-595, arXiv: math.AG/1008.4205v1.

\bibitem{brgh} J. Bryan and A. Gholampour: \emph{The quantum McKay correspondence for polyhedral singularities}, Inventiones Mathematicae, 178(3):655--681, 2009, arXiv: math.AG/0803.3766.

\bibitem{brgr} J. Bryan; T. Graber: \emph{The crepant resolution conjecture}, Algebraic geometry—Seattle 2005. Part 1, 23–42, Proc. Sympos. Pure Math., 80, Part 1, Amer. Math. Soc., Providence, RI, 2009

\bibitem{cale} J. Calabrese: \emph{Donaldson-Thomas invariants and flops}, preprint arXiv: math.AG/1111.1670v1. 

\bibitem{cal2} J. Calabrese: \emph{On the crepant resolution conjecture for Donaldson-Thomas invariants}, preprint arXiv: math.AG/1206.6524v1.

\bibitem{chts} J. C. Chen; H. H. Tseng: \emph{A note on derived McKay correspondence}, Math. Res. Lett. 15 (2008), no. 3, 435–445.

\bibitem{grha} P. Griffiths and J. Harris: \emph{ Principles of algebraic geometry}, Pure and Applied Mathematics, Wiley-Interscience [John Wiley \& Sons], New York, 1978.

\bibitem{gro} A. Grothendieck: \emph{Techniques de construction et theoremes d'existence en geometrie algebrique IV: Les schemas de Hilbert}, Seminaire Bourbaki 1960/61, no 221.

\bibitem{hrs} D. Happel, I. Reiten, S. Smal\o: \emph{Tilting in abelian categories and quasitilted algebras},
Mem. Amer. Math. Soc. 120 (1996), no. 575, viii+ 88 pp. 

\bibitem{hana} G. Harder, and M.S. Narasimhan: \emph{On the cohomology groups of moduli spaces of vector bundles on curves,} Math. Ann. 212 (1975), 215--248.

\bibitem{hule} D. Huybrechts, and M. Lehn: \emph{The geometry of moduli spaces of sheaves}, Aspects of Mathematics, Viewveg 1997. 

\bibitem{joy1} D. Joyce: \emph{Configurations in abelian categories. I. Basic properties and moduli stacks}, Advances in Mathematics 203 (2006), 194-255.

\bibitem{joy2} D. Joyce: \emph{Configurations in abelian categories. II. Ringel-Hall algebras}, Adv. Math. 210 (2007), no. 2, 635–706.

\bibitem{joy3} D. Joyce: \emph{Configurations in abelian categories. III. Stability conditions and identities,} Adv. Math. 215 (2007), no. 1, 153–219. 

\bibitem{joy4} D. Joyce: \emph{Configurations in abelian categories. IV. Invariants and changing stability conditions}, Advances in Mathematics 217 (2008), 125-204.

\bibitem{joy5} D. Joyce: \emph{Motivic invariants of Artin stacks and `stack functions'}, Q. J. Math. 58 (2007), no. 3, 345–392. 

\bibitem{joso} D. Joyce and Y. Song: \emph{A theory of generalized Donaldson-Thomas invariants}, Memoirs of the AMS, 2011,  arXiv: math.AG/0810.5645. 

\bibitem{kle} S. Kleiman: \emph{Les Theoremes de finitude pour le foncteur de Picard}, SGA de Bois Maris 1966/67, exp XIII.

\bibitem{koso} M. Kontsevich, and Y. Soibelman: \emph{Stability structures, motivic Donaldson-Thomas invariants and cluster transformations}, preprint arXiv: math.AG/0811.2435v1. 

\bibitem{lep} J. Le Potier: \emph{Syst\`{e}mes coh\'{e}rents et structures de niveau}, Ast\'{e}risque 214, Soci\'{e}t\'{e} Math\'{e}matique de France, 1993. 

\bibitem{lepa} M. Levine, and R. Pandharipande: \emph{Algebraic cobordism revisited,} arXiv: math.AG/0605196.

\bibitem{lij}  J. Li: \emph{Zero dimensional Donaldson-Thomas invariants of threefolds}, Geom. Topol. 10 (2006) 2117--2171. 

\bibitem{mnop} D. Maulik, N. Nekrasov, A. Okounkov, and R. Pandharipande: \emph{Gromov-Witten theory and Donaldson-Thomas theory}, I. Compos. Math. 142, no. 5, 1263--1285. 

\bibitem{mit} B. Mitchell, \emph{The full imbedding theorem,} Amer. J. Math., 86 (1964) 619–637. 

\bibitem{mum} D. Mumford: \emph{Lectures on curves on an algebraic surface}, Annals Math. Studies 59, Univ Press, Princeton 1966.

\bibitem {oss} C. Okonek, M. Schneider, and H. Spindler: \emph{Vector bundles on complex projective spaces,} Birkhauser, 1980.

\bibitem{pt1} R. Pandharipande, and R. P. Thomas: \emph{Curve counting via stable pairs in the derived category}, Invent. Math. 178 (2009), no. 2, 407--447. 

\bibitem{pt2} R. Pandharipande, and R. P. Thomas: \emph{The 3-fold vertex via stable pairs}, arXiv: math.AG/0709.3823.

\bibitem{pt3} R. Pandharipande, and R. P. Thomas: \emph{Stable pairs and BPS invariants}, arXiv: math.AG/0711.3899.

\bibitem{rei} M. Reineke: \emph{Torsion pairs induced from Harder-Narasimhan filtration}, www.math.uni-bielefeld.de/~sek/select/fahrtorsion.pdf, November, 2007

\bibitem{tho} R. P. Thomas: \emph{A holomorphic Casson invariant for Calabi-Yau 3-folds, and bundles on K3 fibrations,} Journal of Differential Geometry 54, no. 2, 367--438.

\bibitem{to} Y. Toda: \emph{Limit stable objects on Calabi-Yau 3-folds}, Duke Math. J. 149 (2009), no. 1, 157–208.

\bibitem{tod} Y. Toda: \emph{Curve counting theories via stable objects II. DT/ncDT/flop formula}, arxiv: math.AG/0909.5129.

\bibitem{vaf} C. Vafa: \emph{String vacua and orbifoldized LG models} Modern Phys. Lett. A, 4(12):1169--1185, 1989.

\bibitem{van} M. van den Bergh: \emph{Three-dimensional flops and noncommutative rings,} Duke Math. J. 122 (2004), no. 3, 423--455. 

\bibitem{vie}E. Viewhweg: \emph{Rational singularities of higher dimensional schemes}, Proc. Amer. Math. Soc. 63 (1977), no. 1, 6–8. 

\bibitem{weib} C. A. Weibel: \emph{ An introduction to homological algebra,} Cambridge Studies in Advanced Mathematics, 38. Cambridge University Press, Cambridge, 1994. xiv+450 pp.

\bibitem{zas} E. Zaslow: \emph{Topological orbifold models and quantum cohomology rings}, Comm. Math. Phys., 156(2):301--331, 1993.

\end{thebibliography}
\end{document}